\documentclass[a4paper,12pt]{amsart}
\usepackage{amsmath,amsthm}
\usepackage{amsfonts,amsmath,amsthm, amssymb, mathrsfs}
\usepackage{float}
\usepackage{euscript,eufrak,verbatim}
\usepackage{graphicx}
\usepackage[usenames]{color}
\usepackage[colorlinks,linkcolor=red,anchorcolor=blue,citecolor=blue]{hyperref}
\usepackage{amsmath}
\usepackage{amsthm}
\usepackage[all]{xy}
\usepackage{graphicx} 
\usepackage{amssymb} 
\usepackage{enumerate}
\usepackage{enumitem}

\newtheorem{thm}{Theorem}[section]
\newtheorem{prop}[thm]{Proposition}
\newtheorem{df}[thm]{Definition}
\newtheorem{lem}[thm]{Lemma}

\makeatletter
\@addtoreset{equation}{section}
\numberwithin{equation}{section}

\title{Packing topological pressure for amenable group actions}

\author{ Ziqing Ding, Ercai Chen* and Xiaoyao Zhou
}
\address
{1.School of Mathematical Sciences and Institute of Mathematics,  Key Laboratory of NSLSCS, Ministry of Education, Nanjing Normal University, Nanjing 210023, Jiangsu, P.R.China}
\email{mathdingziqing@126.com}
\email{ecchen@njnu.edu.cn}
\email{zhouxiaoyaodeyouxian@126.com}
\date{}

\begin{document}

\maketitle

\renewcommand{\thefootnote}{}
\footnote{*corresponding author}

\begin{abstract}
	In this paper, we first prove the variational principle for amenable packing topological pressure. Then we obtain an inequality concerning amenable packing pressure for factor maps. Finally, we show that the equality about packing topological pressure of the set of generic points when the system satisfies the almost specification property, or $\mu$ is ergodic.
\end{abstract}

\section{Introduction}
Kolomogorov was the first to introduce measure-theoretic entropy in \cite{re20}, and later Adler, Konheim and McAndrew introduced topological entropy in \cite{re21}.
The variational principle establishes a relationship between measure-theoretic entropy and topological entropy in \cite{re16}.
In 1973, Bowen introduced a new definition of topological entropy \cite{re8} on a noncompact set, which was known as Bowen topological entropy or dimensional entropy. Then Pesin and Pitskel \cite{re9} , which we called Pesin-Pitskel pressure. The Hausdorff dimension served as their inspiration. It makes sense to think about similar ideas in dynamical systems with different types of dimensions. For example, the packing dimension can help us come up with a new topological entropy and variational principle for sets that are not compact.
Feng and Huang obtained two variational principles for Bowen entropy and packing entropy in \cite{re10}. It's a cornerstone for the rest of our research. Then Tang, Cheng and Zhao extended the Bowen entropy variational principle to Pensin-Pitskel pressure in \cite{re25}, and Zhong and Chen
extended the packing entropy variational principle to packing pressure in \cite{re2}. Zheng and Chen extended the Bowen entropy variational principle to amenable group actions in \cite{re5} and Dou, Zheng and Zhou extended the packing entropy variational principle to amenable group actions in \cite{re1}, where $G$ is a topological group acting continuously on $X$ instead of $Z.$
For $Z$ action, $T$-invariant Borel probability measure are always exicting. But there are some groups for which there is no invariant probability measure on $X$ with $G$-actions.
In the case where $G$ is amenable, there exists a $G$-invariant Borel probability measure.
In contrast to $Z$ action, a general countable amenable group may have a very complicated structure, which makes
it harder to study. 
Huang, Li and Zhou obtained the variational principle of Pesin-Pitskel pressure on amenable group actions in \cite{re26}.
In this paper we focus on packing topological pressure within the framework of countable discrete amenable group actions.

The first section presents the amenable group and the three theorems that require proof in this paper. The second section mainly introduces packing topological pressure and measure pressure on amenable group actions. Sections 3, 4 and 5 are mainly the proofs of the above three theorems. 
\begin{thm}\label{thm 1.1}
Let $(X,G)$ be a G-action topological dynamical system and G a countable infinite discrete amenable group. Let $\{F_n\}_{n=1}^{\infty}$ be a sequence of finite subsets in G satisfying $\lim_{n \to \infty}\frac{|F_n|}{\log n}=\infty.$ Then for any non-empty analytic subset Z of X and $f \in C(X,\mathbb{R})$.
\\If $P^P(Z,\{F_n\}_{n=1}^{\infty},f)>{\left\lVert f\right\rVert}_\infty.$ 
Then 
\begin{align*}
	P^P(Z,\{F_n\}_{n=1}^{\infty},f)&=\sup\{\overline P_\mu(\{F_n\}_{n=1}^{\infty},f):\mu \in \mathcal{M} (X),\mu(Z)=1\}\\
&=\sup\{P_\mu^P(\{F_n\}_{n=1}^{\infty},f):\mu \in \mathcal{M} (X),\mu(Z)=1\}\\
&=\sup\{P_\mu^{KP}(\{F_n\}_{n=1}^{\infty},f):\mu \in \mathcal{M} (X),\mu(Z)=1\},
\end{align*}
where $P^P(Z,\{F_n\}_{n=1}^{\infty},f),$ $\overline P_\mu(\{F_n\}_{n=1}^{\infty},f),$ $P_\mu^P(\{F_n\}_{n=1}^{\infty},f),$ and $P_\mu^{KP}(\{F_n\}_{n=1}^{\infty},f) $ are defined in sections 2.
\end{thm}
Dou, Zheng and Zhou \cite{re1} prove that the first of these equations is in the case of $f=0.$ For the case $G=Z$, the above equation has been proved by Zhong and Chen in \cite{re2}.

Let $(X,G)$ and $(Y,G)$ be two $G$-action topological dynamical systems. A continuous map $\pi:(X,G) \to (Y,G)$ is called a homomorphism or factor map from $(X,G)$ to $(Y,G)$
if it is onto and $\pi \circ g=g \circ \pi,$ for all $g \in G.$
\begin{thm}\label{thm 1.2}
	Let G be a countable infinite discrete amenable group. Let $\pi:(X,G)\to (Y,G)$ be a factor map. Let $\{F_n\}_{n=1}^{\infty}$ be any tempered F$\o$lner sequence  in $G $ satisfying $\lim_{n \to \infty}\frac{|F_n|}{\log n}=\infty.$ Then for any non-empty subset $E$ of $X$ and $f \in C(X,\mathbb{R})$.
	Then 
	\begin{align*}
	P^P(\pi(E),\{F_n\}_{n=1}^{\infty},f) &\leq P^P(E,\{F_n\}_{n=1}^{\infty},f\circ \pi) \\
	&\leq P^P(\pi(E),\{F_n\}_{n=1}^{\infty},f)+\sup_{y \in Y}h_{top}^{UC}(\pi^{-1}(y),\{F_n\}_{n=1}^{\infty}).
    \end{align*}
\end{thm}
We remark here that for $Z$-actions, Bowen \cite{re17}  proposed the inequality. Then Fang, Huang, Li and Zhang \cite{re27} proved it when entropy is Bowen entropy. Oprocha and Zhang \cite{re31} also proved the above by using open covers.
Li, Chen and Zhou \cite{re28} extended it to Pensin-Pitskel topological pressure.
C. Zhao, E. Chen, X. Hong and X. Zhou proved the inequalities to packing topological pressure in \cite{re6}. And for
amenable group actions, the inequalities of packing entropy were proved by D. Dou, D. Zheng and X. Zhou in \cite{re1}.  
\begin{thm}\label{thm 1.3}
	Let $(X,G)$ be a $G$-action topological dynamical system and $G$ a  countable infinite discrete amenable group. Let $\mu \in \mathcal{M} (X,G)$ and $\{F_n\}_{n=1}^{\infty}$ be a  F$\o$lner sequence in $G$ satisfying $\lim_{n \to \infty}\frac{|F_n|}{\log n}=\infty.$ Then for any  $f \in C(X,\mathbb{R})$.\\
(1)If $(X,G)$ satisfies almost specification property, then 
$$P^P(G_\mu,\{F_n\}_{n=1}^{\infty},f)=h_\mu(X,G)+\int f \,d\mu.$$ 
(2)If $\mu$ is ergodic and $\{F_n\}$ is tempered, then 
$$P^P(G_\mu,\{F_n\}_{n=1}^{\infty},f)=h_\mu(X,G)+\int f \,d\mu.$$ 
\end{thm}
Bowen proved the equality for the case $G=Z$ when $\mu$ is ergodic in \cite{re8}. Pfister and Sullivan extended it when the system satisfies the so-called $g$-almost product property in \cite{re19}.
Zheng and Chen proved the equality of Bowen entropy for amenable group actions in \cite{re4}.
Zhang proved it in \cite{re3} and come to the conclusion that weak specification implies almost specification for amenable systems.
Then, Dou, Zheng and Zhou showed that the packing topological entropy of the set of generic
points for any invariant Borel probability measure $\mu$ coincides with the metric entropy if
either $\mu$ is ergodic or the system satisfies a kind of specification property.

\section{Preliminaries}
In this section, we will introduce amenable group actions and some properties of    topological pressure and measure-theoretic pressure.
\subsection{Amenable group actions}
Let $(X, G)$ be a $G$-action topological dynamical system, where $X$ is a compact metric space with metric $d$ and $G$ is
a topological group acting continuously on $X$. Throughout this paper we assume that $G$
is a countable infinite discrete amenable group. Recall that a
countable discrete group $G$ is amenable if there is a sequence of non-empty finite subsets
$\{F_n\}_{n=1}^{\infty}$ of $G$ which are asymptotically invariant, that is for all $g \in G,$
$$\lim_{n \to \infty} \frac{|F_n \bigtriangleup gF_n|}{|F_n|}=0,$$
which is called a F$\o$lner sequence. Since $G$ is infinite, the sequence $|F_n|$ tends to infinity. Without loss of generality we can assume that $|F_n|$ increases when $n$ increases. One can  refer to \cite{re14,re15}.
A F$\o$lner sequence $\{F_n\}_{n=1}^{\infty}$ is said to be tempered if there exists a constant $C>0$ which is independent of $n$ such that 
$$|\bigcup_{k<n}F_k^{-1}F_n| \leq C|F_n|.$$
Let $\mathcal{M} (X),$ $\mathcal{M} (X, G)$ and $\mathcal{E}(X, G)$ be
the collection of all, $G$-invariant and ergodic $G$-invariant Borel probability measures on $X$,
respectively. Since $G$ is amenable, $\mathcal{M}(X, G)$ and $\mathcal{E}(X, G)$ are both non-empty. For
$\mu \in \mathcal{M}(X, G),$ let $h_\mu(X,G)$ denote the measure-theoretic entropy of $(X, G)$ with respect
to $\mu$.
Let $$G_{\mu,\{F_n\}}=\left\{x \in X:\lim_{n \to \infty}\frac{1}{|F_n|}\sum_{g \in F_n}f(gx)=\int_{X}f \,d\mu ,\forall f \in C(X, \mathbb{R})\right\}$$be the set of generic points of $\mu$ with respect to $\{F_n\}_{n=1}^\infty$.
For simplicity, we can write $G_{\mu,\{F_n\}}$ as $G_{\mu}$. But note that for different F$\o$lner sequence, the corresponding $G_{\mu}$  may not coincide. 

The system $(X,G)$ is said to have the \textit{almost specification  property}  if there exists a non-decreasing function $g:(0,1)\to(0,1)$ with $\lim_{r \to 0}g(r)=0$ and
a map $m:(0,1) \to \mathbb{F}(G) \times (0,1) $ where $\mathbb{F}(G)$ denotes the collection of all finite subsets of $G,$
such that for any $k \in \mathbb{N},$ any $\epsilon_1, \epsilon_2,...,\epsilon_k \in(0,1)$ and any $x_1,x_2,...,x_k \in X,$ if $F_i$ is $m(\epsilon_i)$-invariant, $i=1,2,...,k,$
and $\{F_i\}_{i=1}^k$ are pairwise disjoint, then 
$$\bigcap_{1 \leq i \leq k}B(g;F_i,x_i,\epsilon_i) \neq \emptyset,$$ 
where $B(g;F,x,\epsilon):=\{y \in X:|\{h \in F:d(hx,hy) > \epsilon\}| \leq g(\epsilon)|F|\}$ is the Bowen ball allowing a mistake function with density $g(\epsilon)$ and $F$ is $(K,\delta)$-invariant if $ \frac{\lvert \{g\in G:Kg \cap F \neq \emptyset,Kg \cap (G-F)  \neq \emptyset\} \rvert}{\lvert F \rvert} < \delta.$
 
\subsection{Topological pressure of subsets}
Given $F\in \mathbb{F}(G)$, $x,y \in X$, the Bowen metric $d_F$  on $X$ is defined by  $d_F(x,y):=\max_{g \in F} \limits d(gx,gy).$ Then  \emph{Bowen open  ball } of radius $\epsilon$  in the metric $d_F$ around $x$ is   given by 
$$B_F(x,\epsilon)=\{y\in X: d_F(x,y)<\epsilon\},$$
and \emph{Bowen closed  ball }is given by
$$\overline B_F(x,\epsilon)=\{y\in X: d_F(x,y)\leq\epsilon\}.$$
Given $f\in C(X,\mathbb{R})$, where $C(X,\mathbb{R})$ denotes the set of all continuous functions, define
\begin{align*}
f_F(x)&=\sum_{g \in F}f(gx),\\
f_F(x,\epsilon)&=\sup_{y \in B_F(x,\epsilon)}f_F(y),\\
\overline f_F(x,\epsilon)&=\sup_{y \in \overline B_F(x,\epsilon)}f_F(y).
\end{align*}
Let $Z\subset X$ be a non-empty subset, $\epsilon>0$, $f \in C(X,\mathbb{R})$, $N\in \mathbb{N}$, $s \in \mathbb{R}$ and $\{F_n\}_{n=1}^{\infty}$ be a sequence of finite subsets in $G$ satisfying $|F_n| \to \infty(n \to \infty).$
  Put
$$M^P(N,s,\epsilon,Z,\{F_n\}_{n=1}^{\infty},f)=\sup\{\sum_i\limits  e^{-s\lvert F_{n_i} \rvert +f_{F_{n_i}}(x_i)}\},$$
where the supremum  is taken over all  finite or countable disjoint \\$\{\overline B_{F_{n_i}}(x_i,\epsilon)\}_{i\in I}$ with $n_i \geq N,~x_i \in Z.$  Since $M^P(N,s,\epsilon,Z,\{F_n\}_{n=1}^{\infty},f)$ is decreasing when $N$ increases,
 the following limit exists.
Set
$$M^P(s,\epsilon,Z,\{F_n\}_{n=1}^{\infty},f)=\lim_{N\to \infty}M^P(N,s,\epsilon,Z,\{F_n\}_{n=1}^{\infty},f).$$
Put
\begin{align*}
	M^\mathcal{P} (s,\epsilon,Z,\{F_n\}_{n=1}^{\infty},f)&=\inf \{{\sum_{i=1}^{\infty}M^P(s,\epsilon,Z_i,\{F_n\}_{n=1}^{\infty},f):Z \subset \cup_{i=1}^{\infty} Z_i}\},\\
P^P(\epsilon,Z,\{F_n\}_{n=1}^{\infty},f)&=\sup{\{s:M^\mathcal{P} (s,\epsilon,Z,\{F_n\}_{n=1}^{\infty},f)=\infty\}},\\
&=\inf{\{s:M^\mathcal{P} (s,\epsilon,Z,\{F_n\}_{n=1}^{\infty},f)=0\}},\\
P^P(Z,\{F_n\}_{n=1}^{\infty},f)&=\lim_{\epsilon \to 0}P^P(\epsilon,Z,\{F_n\}_{n=1}^{\infty},f).
\end{align*}
Since $P^P(\epsilon,Z,\{F_n\}_{n=1}^{\infty},f)$ is increasing when $\epsilon$ decreases, the above limit exists. Then we call $P^P(Z,\{F_n\}_{n=1}^{\infty},f)$ {\it packing topological pressure} of the set $Z$ along $\{F_n\}_{n=1}^\infty$  with respect to $f.$
When $f=0$,  $P^P(Z,\{F_n\}_{n=1}^{\infty},0)$ is packing topological entropy $h_{top}^{P}(Z,\{F_n\}_{n=1}^{\infty})$, which is defined by Dou, Zheng and Zhou \cite{re1}.

Let  $Z\subset X$ be a non-empty subset, $\epsilon>0$, $f \in C(X,\mathbb{R})$, $N\in \mathbb{N}$ , $s \in \mathbb{R}$ and $\{F_n\}_{n=1}^{\infty}$ be a sequence of finite subsets in $G$ satisfying $|F_n| \to \infty(n \to \infty).$  
Put
$$M^B(N,s,\epsilon,Z,\{F_n\}_{n=1}^{\infty},f)=\inf\{\sum_i\limits  e^{-s\lvert F_{n_i} \rvert +f_{F_{n_i}}(x_i)}\},$$
where the infimum  is taken over all  finite or countable cover $\{B_{F_{n_i}}(x_i,\epsilon)\}_{i\in I}$ of $Z$ with $n_i \geq N.$
Since  $M^B(N,s,\epsilon,Z,\{F_n\},f)$ is non-decreasing when $N$ increases, the following limit exists.
\begin{align*}
	M^B(s,\epsilon,Z,\{F_n\}_{n=1}^{\infty},f)&=\lim_{N\to \infty}M^B(N,s,\epsilon,Z,\{F_n\}_{n=1}^{\infty},f),\\
P^B(\epsilon,Z,\{F_n\}_{n=1}^{\infty},f)&=\sup{\{s:M^B (s,\epsilon,Z,\{F_n\}_{n=1}^{\infty},f)=\infty\}},\\
&=\inf{\{s:M^B (s,\epsilon,Z,\{F_n\}_{n=1}^{\infty},f)=0\}},\\
P^B(Z,\{F_n\}_{n=1}^{\infty},f)&=\lim_{\epsilon \to 0}P^B(\epsilon,Z,\{F_n\}_{n=1}^{\infty},f).
\end{align*}
 $P^B(Z,\{F_n\}_{n=1}^{\infty},f)$ is called {\it Pesin-Pitskel topological pressure} of the set  $Z$ along  $\{F_n\}_{n=1}^{\infty}$ with respect to $f.$

Let  $Z\subset X$ be a non-empty subset, $F \in \mathbb{F}(G)$, $E\subset X$ and $E$ is  an  {\it $(F,\epsilon)$-spanning set of $Z$} if  for any $x \in Z$, there  exists  $y\in E$ such that $d_F(x,y)< \epsilon.$  
A set $E\subset Z$ is an  {\it  $(F,\epsilon)$-separated set of $Z$} if  for any $x, y \in E$ with $x\not=y$, one has  $d_F(x,y)> \epsilon.$  Let $\{F_n\}_{n=1}^\infty\subset\mathbb{F}(G)$  satisfying $|F_n| \to \infty(n \to \infty).$
Put
$$Q(Z,\epsilon,{F_n},f)=\inf{\{\sum_{x \in E}e^{f_{F_n}(x)}:\emph{E is an $(F_n,\epsilon)$-spanning set of $Z$}\}},$$
$$P(Z,\epsilon,{F_n},f)=\sup{\{\sum_{x \in E}e^{f_{F_n}(x)}:\emph{E is an $(F_n,\epsilon)$-separated set of $Z$}\}}.$$
Let $\delta >0$ and choose $\epsilon >0$ so that $d(x,y)<\frac{\epsilon}{2}$ implies $|f(x)-f(y)| < \delta.$ 
It is easy to get
\begin{align}\label{a}
Q(Z,\epsilon,{F_n},f) \leq P(Z,\epsilon,{F_n},f) \leq e^{|F_n|\delta}Q(Z,\frac{\epsilon}{2},{F_n},f). 
\end{align}
Set
\begin{align*}
	P^{UC}(Z,\epsilon,\{F_n\}_{n=1}^{\infty},f)&=\limsup_{n \to \infty}  \frac{1}{\lvert F_n \rvert} \log P(Z,\epsilon,F_n,f),\\
	P^{UC'}(Z,\epsilon,\{F_n\}_{n=1}^{\infty},f)&=\limsup_{n \to \infty}  \frac{1}{\lvert F_n \rvert} \log Q(Z,\epsilon,F_n,f).
\end{align*}
By \ref{a}, we have
\begin{align*}
	P^{UC'}(Z,\epsilon,\{F_n\}_{n=1}^{\infty},f) &\leq P^{UC}(Z,\epsilon,\{F_n\}_{n=1}^{\infty},f) \\&\leq P^{UC'}(Z,\frac{\epsilon}{2},\{F_n\}_{n=1}^{\infty},f)+\delta.
\end{align*}
Letting $\epsilon\to0, $ we can get
\begin{align*}
P^{UC}(Z,\{F_n\}_{n=1}^{\infty},f)&=\lim_{\epsilon \to 0}P^{UC}(Z,\epsilon,\{F_n\}_{n=1}^{\infty},f)\\&=\lim_{\epsilon \to 0}P^{UC'}(Z,\epsilon,\{F_n\}_{n=1}^{\infty},f).
\end{align*} 
$P^{UC}(Z,\{F_n\}_{n=1}^{\infty},f)$  is called {\it upper capacity topological pressure} of  the set $Z$ along $\{F_n\}_{n=1}^{\infty}$  with respect to $f.$
\subsection{Properties of amenable topological pressure}
The following propositions are some basic properties of topological pressure. Below we will prove that Zhang's definition of {\it Pesin-Pitskel topological pressure} \cite{re3} is consistent with our definition above, which will help in proving Theorem \ref{thm 1.3}.
\begin{prop}\label{prop 2.1}
	If we define  $P^{P'}(\epsilon,Z,\{F_n\}_{n=1}^{\infty},f)$ and $P^{B'}(\epsilon,Z,\{F_n\}_{n=1}^{\infty},f),$ by
replacing the above $f_{F_{n_i}}(x_i)$  by $f_{F_{n_i}}(x_i,\epsilon)$ and $\overline f_{F_{n_i}}(x_i,\epsilon)$, then 
$$P^{P}(Z,\{F_n\}_{n=1}^{\infty},f)=\lim_{\epsilon \to 0}P^{P'}(\epsilon,Z,\{F_n\}_{n=1}^{\infty},f),$$
$$P^{B}(Z,\{F_n\}_{n=1}^{\infty},f)=\lim_{\epsilon \to 0}P^{B'}(\epsilon,Z,\{F_n\}_{n=1}^{\infty},f).$$
\end{prop}

\begin{proof}
Fix $\epsilon>0$. It is clear that $$P^{B}(\epsilon,Z,\{F_n\}_{n=1}^{\infty},f) \leq P^{B'}(\epsilon,Z,\{F_n\}_{n=1}^{\infty},f).$$
Let $$var(f,\epsilon)=\sup \{\lvert f(x)-f(y)\rvert\ :d(x,y) \leq \epsilon\}.$$
This shows that $$f_{F_{n}}(x,\epsilon)-\lvert F_{n}\rvert var(f,\epsilon) \leq f_{F_{n}}(x).$$

\begin{align*}
M^B(N,s,\epsilon,Z,\{F_n\}_{n=1}^{\infty},f)&=\inf\{\sum_i\limits  e^{-s\lvert F_{n_i} \rvert +f_{F_{n_i}}(x_i)} \}\\
&\geq\inf\{\sum_i\limits  e^{-s\lvert F_{n_i} \rvert +f_{F_{n_i}}(x_i,\epsilon)-\lvert F_{n_i} \rvert var(f,\epsilon)} \}\\
&= M^{B'}(N,s+var(f,\epsilon),\epsilon,Z,\{F_n\}_{n=1}^{\infty},f).	
\end{align*}
Letting $N \to \infty$, we have $$M^B(s,\epsilon,Z,\{F_n\}_{n=1}^{\infty},f)\geq M^{B'}(s+var(f,\epsilon),\epsilon,Z,\{F_n\}_{n=1}^{\infty},f),$$
$$P^{B}(\epsilon,Z,\{F_n\}_{n=1}^{\infty},f) \geq P^{B'}(\epsilon,Z,\{F_n\}_{n=1}^{\infty},f)-var(f,\epsilon).$$
Letting $\epsilon \to \infty$, we get $$P^{B}(Z,\{F_n\}_{n=1}^{\infty},f)=\lim_{\epsilon \to 0}P^{B'}(\epsilon,Z,\{F_n\}_{n=1}^{\infty},f).$$ The next equality can be proved similarly.
\end{proof}

The following four properties are obvious and we will omit the proof.
\begin{prop}\label{prop 2.2}
Let $\{F_n\}_{n=1}^{\infty}$ be a F$\o$lner sequence in $G$.\\
(1)If $Z_1 \subset Z_2,$ then for any $\epsilon>0$, $f \in C(X,\mathbb{R})$, $N\in \mathbb{N}$, $s \in \mathbb{R},$ $M^P(N,s,\epsilon,Z_1,\{F_n\}_{n=1}^{\infty},f) \leq M^P(N,s,\epsilon,Z_2,\{F_n\}_{n=1}^{\infty},f),$ and moreover $M^P(s,\epsilon,Z_1,\{F_n\}_{n=1}^{\infty},f) \leq M^P(s,\epsilon,Z_2,\{F_n\}_{n=1}^{\infty},f).$\\
(2)If $Z_1 \subset Z_2,$ then $\mathfrak{P} (Z_1,\{F_n\}_{n=1}^{\infty},f)\leq\mathfrak{P} (Z_2,\{F_n\}_{n=1}^{\infty},f),$  where $\mathfrak{P} \in \{P^B,P^P,P^{UC}\}.$\\
(3)For all $c \in  \mathbb R,$ $\mathfrak{P} (Z,\{F_n\}_{n=1}^{\infty},f+c)=\mathfrak{P} (Z,\{F_n\}_{n=1}^{\infty},f)+c,$  where $\mathfrak{P} \in \{P^B,P^P,P^{UC}\}.$\\
(4)If $Z \subset \cup_{i=1}^{\infty}Z_i$, then for any $\epsilon>0,$ $s \in \mathbb R, f\in C(X,\mathbb R),$ 
$$M^\mathcal{P} (s,\epsilon,\cup_{i} Z_i,\{F_n\}_{n=1}^{\infty},f) \leq \sum_{i=1}^{\infty}M^\mathcal{P} (s,\epsilon,Z_i,\{F_n\}_{n=1}^{\infty},f),$$
$$P^P(\epsilon,Z,\{F_n\}_{n=1}^{\infty},f) \leq \sup_{i}P^P(\epsilon,Z_i,\{F_n\}_{n=1}^{\infty},f),$$ 
$$P^P(Z,\{F_n\}_{n=1}^{\infty},f) \leq \sup_{i}P^P(Z_i,\{F_n\}_{n=1}^{\infty},f).$$ 
\end{prop}

Next we will compare {\it packing topological pressure} with {\it Pesin-Pitskel topological pressure}. Dou, Zheng and Zhou proved that {\it Bowen topological entropy} was smaller than {\it packing topological entropy} in \cite{re1} for amenable group actions. We now extend it to the pressure.
\begin{prop}\label{prop 2.3}
	Let $\{F_n\}_{n=1}^{\infty}$ be a F$\o$lner sequence in $G$.  For any $Z \subset X,$
	$$P^{B}(Z,\{F_n\}_{n=1}^{\infty},f) \leq P^{P}(Z,\{F_n\}_{n=1}^{\infty},f).$$
\end{prop}
The proof is similar to the proof of Proposition 2.4(4) in \cite{re2}. We give the proof here for completeness.
\begin{proof}
	Suppose that $P^{B}(Z,\{F_n\}_{n=1}^{\infty},f)>s>-\infty$. For any $\epsilon>0$ $\empty$ and $\empty$ $n\in \mathbb{N}$, let 
	$$\mathcal{F} _{F_n,\epsilon}=\{\mathcal{F}:\mathcal{F}=\{\overline B_{F_n}(x_i,\epsilon)\}disjoint,x_i \in Z \}.$$
	Take $\mathcal{F}({F_n},\epsilon,Z) \in \mathcal{F} _{F_n,\epsilon}$ such that $\lvert \mathcal{F}({F_n},\epsilon,Z) \rvert=\max_{\mathcal{F} \in \mathcal{F} _{F_n,\epsilon}}{\lvert \mathcal{F} \rvert}.$
	For convenience, we denote $\mathcal{F}({F_n},\epsilon,Z)=\{\overline B_{F_n}(x_i,\epsilon):i=1,...,\lvert \mathcal{F}({F_n},\epsilon,Z) \rvert\}.$
	It is easy to check that
	$$Z \subset \bigcup\limits_ {i=1}^{\lvert \mathcal{F}({F_n},\epsilon,Z) \rvert}{B_{F_n}}(x_i,2\epsilon+\delta),\forall \delta>0.$$
	Then for any $s \in \mathbb{R},$

	$$M^B(n,s,2\epsilon+\delta,Z,\{F_n\}_{n=1}^{\infty},f) \leq  e^{-s\lvert{F_n}\rvert}\sum_{i=1}^{\lvert \mathcal{F}({F_n},\epsilon,Z) \rvert}e^{f_{F_n}(x_i)}$$
	$$\leq M^P(n,s,\epsilon,Z,\{F_n\}_{n=1}^{\infty},f).$$

This implies that $$M^B(s,2\epsilon+\delta,Z,\{F_n\}_{n=1}^{\infty},f) \leq M^\mathcal{P} (s,\epsilon,Z,\{F_n\}_{n=1}^{\infty},f).$$
	Since $P^{B}(Z,\{F_n\}_{n=1}^{\infty},f)>s>-\infty$, $M^B(s,2\epsilon+\delta,Z,\{F_n\}_{n=1}^{\infty},f) \geq 1$ when $\epsilon$  and $\delta$ are small enough, furthermore, $M^\mathcal{P} (s,\epsilon,Z,\{F_n\}_{n=1}^{\infty},f) \geq 1.$
	This shows
	$P^P(\epsilon,Z,\{F_n\}_{n=1}^{\infty},f) \geq s.$ Letting $s\to P^{B}(Z,\{F_n\}_{n=1}^{\infty},f), $ we have
	$P^{B}(Z,\{F_n\}_{n=1}^{\infty},f) \leq P^{P}(Z,\{F_n\}_{n=1}^{\infty},f).$
\end{proof}

Naturally we will compare {\it packing topological pressure} with {\it upper capacity topological pressure}. We can see similar conclusions in \cite{re1, re2}.
\begin{prop}\label{prop 2.4}
	Let $\{F_n\}_{n=1}^{\infty}$ be a F$\o$lner sequence in $G$ satisfying $\lim_{n \to \infty}\frac{|F_n|}{\log n}=\infty.$ Then for any subset $Z \subset X,$ $f \in C(X,\mathbb{R})$ and any $\epsilon >0,$
$$P^{P}(\epsilon,Z,\{F_n\}_{n=1}^{\infty},f) \leq P^{UC}(Z,\epsilon,\{F_n\}_{n=1}^{\infty},f).$$
Furthermore,
$$P^{P}(Z,\{F_n\}_{n=1}^{\infty},f) \leq P^{UC}(Z,\{F_n\}_{n=1}^{\infty},f).$$	
\end{prop}
\begin{proof}
	Let $\epsilon>0$ and $-\infty<t<s<P^{P}(\epsilon,Z,\{F_n\}_{n=1}^{\infty},f)$. Then
	$$M^{P} (s,\epsilon,Z,\{F_n\}_{n=1}^{\infty},f) \geq M^\mathcal{P} (s,\epsilon,Z,\{F_n\}_{n=1}^{\infty},f) =\infty.$$
	For any $N \in \mathbb{N}$, there exists a pairwise family $\overline B_{F_{n_i}}(x_i,\epsilon)$ such that $x_i \in Z,$
	$n_i \geq N$, and $\sum_{i}e^{-s\lvert F_{n_i} \rvert +f_{F_{n_i}}(x_i)}>1$.
	For each $k$, let
	$$m_k=\{x_i:n_i=k\}.$$
	Then
	\begin{align}\label{a2}
		\sum_{k=N}^{\infty}\sum_{x \in m_k}e^{f_{F_k}(x)}e^{-\lvert {F_k} \rvert s}>1.
	\end{align}
	Since $\{F_n\}$ satisfies the growth condition $\lim_{n \to \infty}\frac{|F_n|}{\log n}=\infty,$ $\sum_{k=1}^{\infty}e^{\lvert {F_k} \rvert (t-s)}$ converges.
    Let $M=\sum_{k=1}^{\infty}e^{\lvert {F_k} \rvert (t-s)}.$ There must be some $k \geq N$ such that $\sum_{x \in m_k}e^{f_{F_k}(x)}>\frac{1}{M}e^{\lvert {F_k} \rvert t},$
	otherwise the above sum is at most
	$$\sum_{k=1}^{\infty}\frac{1}{M}e^{\lvert {F_k} \rvert t}e^{-\lvert {F_k} \rvert s}=1.$$ This contradicts \ref{a2}.
	So $P(Z,\epsilon,F_k,f) \geq \sum_{x \in m_k}e^{f_{F_k}(x)}>\frac{1}{M}e^{\lvert {F_k} \rvert t}$ and hence
	$$P^{UC}(Z,\epsilon,\{F_n\}_{n=1}^{\infty},f)=\limsup_{k \to \infty}\frac{1}{\lvert {F_k} \rvert} \log P(Z,\epsilon,F_k,f) \geq t.$$
 This gets the first inequality by letting   $t\to P^{P}(\epsilon, Z,\{F_n\}_{n=1}^{\infty},f),$ and  the second inequality by letting $\epsilon \to0 .$  
\end{proof}

The following proposition played a key role in proving Theorem \ref{thm 1.2} and the proof is inspired by \cite[Proposition 2.7]{re1}.
\begin{prop}\label{prop 2.5}
Let $\epsilon>0$, $Z\subset X$ and $\{F_n\}_{n=1}^{\infty}\subset\mathbb{F}(G)$  satisfying $\lim_{n \to \infty}\frac{|F_n|}{\log n}=\infty$, $f \in C(X,\mathbb{R}).$
\begin{flushleft}
(1)$P^P(Z,\{F_n\}_{n=1}^{\infty},f) \leq \inf \{\sup_{i \geq 1}P^{UC}(Z_i,\{F_n\}_{n=1}^{\infty},f):Z=\cup_{i=1}^{\infty}Z_i \}.$\\
(2)For any $ \epsilon >0, \delta >0,$ there exists a cover $ \cup_{i=1}^{\infty}Z_i=Z ,$ such that
$$P^P(Z,\epsilon,\{F_n\}_{n=1}^{\infty},f) +\delta \geq \sup_{i \geq 1}P^{UC}(Z_i,3\epsilon,\{F_n\}_{n=1}^{\infty},f).$$
\end{flushleft}
\end{prop}
\begin{proof}
$(1)$For $Z=\bigcup_{i=1}^{\infty}Z_i,$ we have
\begin{align*}
	&P^P(Z,\epsilon,\{F_n\}_{n=1}^{\infty},f)\\ \leq& \sup_{i \geq 1}P^P(Z_i,\epsilon,\{F_n\}_{n=1}^{\infty},f) (\text{by Proposition } \ref{prop 2.2})\\ \leq &\sup_{i \geq 1}P^{UC}(Z_i,\epsilon,\{F_n\}_{n=1}^{\infty},f)(\text{by Proposition } \ref{prop 2.4})
	\\ \leq &\sup_{i \geq 1}P^{UC}(Z_i,\{F_n\}_{n=1}^{\infty},f).
\end{align*}

 So
$$P^P(Z,\{F_n\}_{n=1}^{\infty},f) \leq \inf \{\sup_{i \geq 1}P^{UC}(Z_i,\{F_n\}_{n=1}^{\infty},f):Z=\cup_{i=1}^{\infty}Z_i \}.$$
$(2)$ Assume that $P^P(Z,\epsilon,\{F_n\}_{n=1}^{\infty},f)<\infty$. Given $ \delta >0$ and $s:=\\P^P(Z,\epsilon,\{F_n\}_{n=1}^{\infty},f)+\delta.$ This means 
$ M^\mathcal{P} (s,\epsilon,Z,\{F_n\}_{n=1}^{\infty},f)=0.$
There exists a cover of $Z \subset \bigcup_{i=1}^{\infty }Z_i'$ such that
$$\sum_{i}M^{P} (s,\epsilon,Z_i',\{F_n\}_{n=1}^{\infty},f)<1 .$$
Set $Z_i=Z_i' \cap Z.$ We have
$$\sum_{i}M^{P} (s,\epsilon,Z_i,\{F_n\}_{n=1}^{\infty},f)<1 .$$
For each $Z_i$, when $N$ is large enough, we have
$$M^{P} (N,s,\epsilon,Z_i,\{F_n\}_{n=1}^{\infty},f)<1 .$$
Let $E$ be an  $(F_N,3\epsilon)$-separated set of $Z_i.$
Since $\{\overline B_{F_N}(x_i,\epsilon):x_i\in  E\}$  is pairwise disjoint, we have
\begin{align*}
 \sum_{x_i \in E}\limits  e^{-s\lvert F_N \rvert +f_{F_N}(x_i)} &\leq M^{P} (N,s,\epsilon,Z_i,\{F_n\}_{n=1}^{\infty},f)<1,\\
  \sum_{x_i \in E}\limits  e^{f_{F_N}(x_i)} &\leq  e^{s\lvert F_N \rvert },\\
  P^{UC}(Z_i,3\epsilon,\{F_n\}_{n=1}^{\infty},f) &\leq s. 
\end{align*}
Furthermore,
$$P^P(Z,\epsilon,\{F_n\}_{n=1}^{\infty},f) +\delta \geq \sup_{i \geq 1}P^{UC}(Z_i,3\epsilon,\{F_n\}_{n=1}^{\infty},f).$$
\end{proof}

\subsection{Measure-theoretic pressure}
Let $f \in C(X,\mathbb{R})$ and $\mu \in \mathcal{M} (X)$. Let $\{F_n\}_{n=1}^{\infty}$ be a sequence of finite subsets of $G$ with $\lvert F_n\rvert \to \infty.$
The {\it measure-theoretic lower and upper local pressure} of $x \in X$ with respect to $\mu$ and $f$ are defined by
$$\underline P_\mu(x,\{F_n\}_{n=1}^{\infty},f)=\lim_{\epsilon \to 0}\liminf_{n \to \infty}\dfrac{-\log\mu(B_{F_n}(x,\epsilon))+f_{F_n}(x)}{\lvert F_n\rvert} .$$
$$\overline P_\mu(x,\{F_n\}_{n=1}^{\infty},f)=\lim_{\epsilon \to 0}\limsup_{n \to \infty}\dfrac{-\log\mu(B_{F_n}(x,\epsilon))+f_{F_n}(x)}{\lvert F_n\rvert} .$$
\begin{df}
	The measure-theoretic lower and upper local pressure with respect to $\mu$ and $f$ are defined by
	$$\underline P_\mu(\{F_n\}_{n=1}^{\infty},f):=\int \underline P_\mu(x,\{F_n\}_{n=1}^{\infty},f)d\mu(x) .$$
	$$\overline P_\mu(\{F_n\}_{n=1}^{\infty},f):=\int \overline P_\mu(x,\{F_n\}_{n=1}^{\infty},f)d\mu(x)  .$$	
\end{df}

The following two definitions were inspired by \cite{re2} and we will prove in Proposition \ref{prop 2.10} that the two of them are consistent.  For the case $G=Z,$ these was defined in \cite{re2}.
\begin{df}
We call the following quantity
$$P_\mu^P(\{F_n\}_{n=1}^{\infty},f):=\lim_{\epsilon \to 0}\lim_{\delta \to 0}\inf\{P^P(\epsilon,Z,\{F_n\}_{n=1}^{\infty},f):\mu(Z)\geq1-\delta\}$$
packing pressure of $\mu$ along  $\{F_n\}_{n=1}^{\infty}$ with respect to $f.$
\end{df}

\begin{df}
Let $$M_\mu^{\mathcal{P} }(s,\epsilon,\delta,\{F_n\}_{n=1}^{\infty},f)=\inf\{\sum_{i=1}^{\infty}M^P(s,\epsilon,Z_i,\{F_n\}_{n=1}^{\infty},f):\mu(\cup_{i=1}^\infty Z_i)\geq 1-\delta\}$$ and
$$P_\mu^{KP}(\epsilon,\delta,\{F_n\}_{n=1}^{\infty},f)=\sup\{s:M_\mu^{\mathcal{P} }(s,\epsilon,\delta,\{F_n\}_{n=1}^{\infty},f)=+\infty\}.$$ We call the following quantity
$$P_\mu^{KP}(\{F_n\}_{n=1}^{\infty},f)=\lim_{\epsilon \to 0}\lim_{\delta \to 0}P_\mu^{KP}(\epsilon,\delta,\{F_n\}_{n=1}^{\infty},f)$$ packing pressure of $\mu$ in the sense of Katok along $\{F_n\}_{n=1}^{\infty}$ with respect to $f.$
\end{df}

\subsection{Properties of amenable measure-theoretic pressure}
\begin{prop}\label{prop 2.9}
Replacing $f_{F_{n_i}}(x_i)$ above by $f_{F_{n_i}}(x_i,\epsilon)$ and $\overline f_{F_{n_i}}(x_i,\epsilon)$ respectively,
we can define new functions $\mathcal{M}_\mu^\mathcal{P}$ and $P_\mu^{KP'}(\epsilon,\delta,\{F_n\}_{n=1}^{\infty},f).$
Then
$$P_\mu^{KP}(\{F_n\}_{n=1}^{\infty},f)=\lim_{\epsilon \to 0}\lim_{\delta \to 0}P_\mu^{KP'}(\epsilon,\delta,\{F_n\}_{n=1}^{\infty},f).$$
\end{prop}
\begin{proof}
The proof is analogous to that of Proposition \ref{prop 2.1}, so we omit it.
\end{proof}
\begin{prop}\label{prop 2.10}
Let $\mu \in \mathcal{M} (X)$, $f \in C(X,\mathbb{R})$ and $\{F_n\}_{n=1}^\infty\subset\mathbb F(G)$ with $\lvert F_n\rvert \to \infty.$
Then
$$P_\mu^{KP}(\{F_n\}_{n=1}^{\infty},f)=P_\mu^{P}(\{F_n\}_{n=1}^{\infty},f).$$
\end{prop}

\begin{proof}
	We first prove that $P_\mu^{KP}(\{F_n\}_{n=1}^{\infty},f) \leq P_\mu^{P}(\{F_n\}_{n=1}^{\infty},f).$
	For any $s < P_\mu^{KP}(\{F_n\}_{n=1}^{\infty},f),$ there exists $\epsilon' >0$ and $\delta'>0$ such that for any $\epsilon \in (0,\epsilon')$ and $\delta \in (0,\delta'),$
	$$P_\mu^{KP}(\epsilon,\delta,\{F_n\}_{n=1}^{\infty},f)>s.$$
	Thus
	$$ M^\mathcal{P} (s,\epsilon,\delta,\{F_n\}_{n=1}^{\infty},f)=\infty.$$
 If $Z\subset \cup_{i=1}^\infty Z_i$ with $\mu(Z) \geq 1-\delta,$   then $\mu(\cup_{i=1}^\infty Z_i) \geq 1-\delta.$
	It follows that
	$$\sum_{i=1}^\infty M^{P} (s,\epsilon,Z_i,\{F_n\}_{n=1}^{\infty},f)=\infty,$$
	which implies that $M^\mathcal{P} (s,\epsilon,Z,\{F_n\}_{n=1}^{\infty},f)=\infty.$
	Hence $P^{P}(\epsilon,Z,\{F_n\}_{n=1}^{\infty},f) \geq s$ and $P_\mu^{P}(\{F_n\}_{n=1}^{\infty},f) \geq s.$ Letting $s\to P_\mu^{KP}(\{F_n\}_{n=1}^{\infty},f), $ this shows that $P_\mu^{KP}(\{F_n\}_{n=1}^{\infty},f) \leq P_\mu^{P}(\{F_n\}_{n=1}^{\infty},f).$
	Next, we shall show the inverse inequality. If $s<P^{P}(\{F_n\}_{n=1}^{\infty},f),$ then there exists $\epsilon'>0$ and $\delta' >0$ such that for any $\epsilon \in (0,\epsilon')$ and $\delta \in (0,\delta'),$
	for any family $\{Z_i\}_{i=1}^\infty$ with $\mu(\cup_{i=1}^\infty Z_i) \geq 1-\delta,$ we have
	$$P^{P}(\epsilon,\cup_{i=1}^\infty Z_i,\{F_n\}_{n=1}^{\infty},f)>s.$$
	This implies that
	$$M^\mathcal{P} (s,\epsilon,\cup_{i=1}^\infty Z_i,\{F_n\}_{n=1}^{\infty},f)=\infty.$$
	Thus
	$$\sum_{i=1}^\infty M^{P} (s,\epsilon,Z_i,\{F_n\}_{n=1}^{\infty},f)=\infty.$$
	Thus
	$$ M^\mathcal{P} (s,\epsilon,\delta,\{F_n\}_{n=1}^{\infty},f)=\infty.$$
	Hence
	$$P_\mu^{KP}(\epsilon,\delta,\{F_n\}_{n=1}^{\infty},f)>s.$$
	Letting $\delta\to0,\epsilon\to0, s\to P^{P}(\{F_n\}_{n=1}^{\infty},f),$ we finish the proof.
\end{proof}

Before proving the next proposition, we need the following classical $5r$-lemma in geometric measure theory(\cite[Theorem 2.1]{re40}).
\begin{lem}\label{lem 2.11}
Let $(X,d)$ be a compact metric space and $\mathcal{B}=\{B(x_i,r_i)\}_{i\in\mathfrak{I} } $be a family of closed (or open) balls in $X.$ Then there exists a finite or countable subfamily
$\mathfrak{B}'=\{B(x_i,r_i)\}_{i\in\mathfrak{I}' }$ of pairwise disjoint balls in $\mathcal{B}$ such that 
$$\bigcup_{B \in \mathfrak{B}}B \subset \bigcup_{i \in \mathfrak{I}' }B(x_i,5r_i).$$
\end{lem}

\begin{prop}\label{prop 2.12}
Let $\mu \in \mathcal{M} (X)$ , $Z \subset X,$ $f \in C(X,\mathbb{R})$ and $\{F_n\}_{n=1}^{\infty}\subset\mathbb F(G)$  satisfying $\lim_{n \to \infty}\frac{|F_n|}{\log n}=\infty.$ 
For $s \in \mathbb{R},$ the following properties hold:\\
(1)If $\overline P_\mu(x,\{F_n\}_{n=1}^{\infty},f) \leq s$ for any  $x \in Z,$ then $P^P(Z,\{F_n\}_{n=1}^{\infty},f) \leq s.$\\
(2)If $\overline P_\mu(x,\{F_n\}_{n=1}^{\infty},f) \geq s$ for any  $x \in Z $  and $\mu(Z)>0,$ then $P^P(Z,\{F_n\}_{n=1}^{\infty},f) \geq s.$	

\end{prop}
\begin{proof}
	We now prove (1). Fix  $\beta >s$ and let 
	$$Z_m=\{x \in Z:\limsup_{n \to \infty}\frac{-\log\mu(B_{F_n}(x,\epsilon))+f_{F_n}(x)}{\lvert F_n \rvert} < \frac{\beta + s}{2},\forall \epsilon \in (0,\frac{1}{m}) \} .$$
	It is easy to check that $Z=\bigcup_{m=1}^{\infty}Z_m.$ 
	Given $m \geq 1,$  $x \in Z_m $  and $\epsilon \in (0,\frac{1}{m}),$ there exists $N \in \mathbb{N}$ such that for any $n \geq N,$
	$$\mu(B_{F_n}(x,\epsilon)) >  e^{-\frac{\beta+s}{2}\lvert F_n \rvert + f_{F_n}(x)}.$$
	Let
	$$Z_{m,N}=\{x \in Z_m:\mu(B_{F_n}(x,\epsilon)) \geq  e^{-\frac{\beta+s}{2}\lvert F_n \rvert + f_{F_n}(x)},\forall n \geq N,\epsilon \in (0,\frac{1}{m})\}.$$
	It is clear that $Z_m=\bigcup_{N=1}^{\infty}Z_{m,N}.$
	Given $\epsilon>0,$ $N\in\mathbb N$ and $L \geq N.$
Let  $\mathcal{F} =\{B_{F_{n_i}}(x_i,\epsilon)\}_{i \in I},$ where $x_i \in Z_{m,N},$ $n_i \geq L $  be a finite or countable disjoint family.
	\begin{align*}
		\sum_{i}e^{-\beta\lvert F_{n_i} \rvert + f_{F_{n_i}}(x_i)}&=\sum_{i}e^{-\lvert F_{n_i} \rvert(\frac{\beta+s}{2}+\frac{\beta-s}{2}) + f_{F_{n_i}}(x_i)}\\
	&\leq e^{-\lvert F_{L} \rvert(\frac{\beta-s}{2})}\sum_{i}e^{-\lvert F_{n_i} \rvert\frac{\beta+s}{2} + f_{F_{n_i}}(x_i)}\\
	&\leq e^{-\lvert F_{L} \rvert(\frac{\beta-s}{2})}\sum_{i}\mu(\overline B_{F_{n_i}}(x_i,\epsilon))\\
	&\leq e^{-\lvert F_{L} \rvert(\frac{\beta-s}{2})}.
	\end{align*}
   It follows that 
	$$M^P(L,\beta,\epsilon,Z_{m,N},\{F_n\}_{n=1}^{\infty},f) \leq e^{-\lvert F_{L} \rvert(\frac{\beta-s}{2})}.$$
	Letting $L\to\infty,$ we have $M^P(\beta,\epsilon,Z_{m,N},\{F_n\}_{n=1}^{\infty},f)=0,$ which implies that $M^\mathcal{P} (\beta,\epsilon,Z_{m},\{F_n\}_{n=1}^{\infty},f)=0.$
	Hence  $P^P(\epsilon,Z_{m},\{F_n\}_{n=1}^{\infty},f) \leq \beta.$
	Letting $\epsilon\to0,$ it follows that $P^P(Z_m,\{F_n\}_{n=1}^{\infty},f) \leq \beta.$
	Combining this and $P^P(Z,\{F_n\}_{n=1}^{\infty},f) \leq \sup_m P^P(Z_m,\{F_n\}_{n=1}^{\infty},f),$ we get $P^P(Z,\{F_n\}_{n=1}^{\infty},f)\leq\beta.$ Letting $\beta\to s,$ we prove (1).\\
Now, we prove (2).
	Fix $\beta <s.$ Let $\delta=\frac{s-\beta}{2}$
and
	$Z_m=\{x \in Z:\limsup_{n \to \infty}\frac{-\log\mu(B_{F_n}(x,\frac{1}{m}))+f_{F_n}(x)}{\lvert F_n \rvert} > \beta + \delta\}.$ Then $Z=\bigcup_{m=1}^{\infty}Z_m.$
	Since $\mu(Z)>0$ and $Z_n\subset Z_{n+1}, n\in\mathbb N,$ there exists $m\in\mathbb N$ such that $\mu(Z_m)>0.$ For any  $\epsilon \in (0,\frac{1}{m})$ and $ x \in Z_m$, we have
	$$ \limsup_{n \to \infty}\frac{-\log\mu(B_{F_n}(x,\epsilon))+f_{F_n}(x)}{\lvert F_{n} \rvert}>\beta +\delta.$$
	Next we claim that $M^\mathcal{P} (s,\frac{\epsilon}{10},Z_m,\{F_n\}_{n=1}^{\infty},f)=\infty,$ which implies that 
	\begin{align*}
		P^P(Z,\{F_n\}_{n=1}^{\infty},f) &\geq P^P(Z_m,\{F_n\}_{n=1}^{\infty},f)\\ &\geq P^P(\frac{\epsilon}{10},Z_m,\{F_n\}_{n=1}^{\infty},f)\\ &\geq s.
	\end{align*}
To this end, it suffices to show that $M^{P} (s,\frac{\epsilon}{10},E,\{F_n\}_{n=1}^{\infty},f)=\infty$ for any Borel subset $E \subset Z_m$ with $\mu (E)>0.$
	In fact, for  $E\subset Z_m$ with $\mu (E)>0,$
	let 
	$E_n=\{x\in E:\mu(B_{F_n}(x,\epsilon))<e^{-\lvert F_{n} \rvert(\beta+\delta)+f_{F_n}(x)}\},~n \in \mathbb{N}.$
	It is clear that $E=\bigcup_{n=N}^{\infty}E_n$ for each $N \in \mathbb{N}. $ Then $\mu(\bigcup_{n=N}^{\infty}E_n)=\mu(E).$
	Hence there exists $n \geq N$ such that
	$$\mu(E_n) \geq \frac{1}{n(n+1)}\mu(E).$$
	Fix such $n$ and let $\mathcal{B} =\{\overline B_{F_n}(x,\frac{\epsilon}{10}):x\in E_n\}.$
	By Lemma \ref{lem 2.11}, there exists a finite or countable pairwise disjoint family $\{\overline B_{F_n}(x_i,\frac{\epsilon}{10})\}_{i\in I}$ such that
	$$E_n \subset \bigcup_{x \in E_n}\overline B_{F_n}(x,\frac{\epsilon}{10}) \subset \bigcup_{i\in I} \overline B_{F_n}(x_i,\frac{\epsilon}{2}) \subset \bigcup_{i\in I} B_{F_n}(x_i,\epsilon).$$
	Hence,
	\begin{align*}
		M^P(N,\beta,\frac{\epsilon}{10},E,\{F_n\}_{n=1}^{\infty},f) &\geq M^P(N,\beta,\frac{\epsilon}{10},E_n,\{F_n\}_{n=1}^{\infty},f)\\ &\geq \sum_{{i\in I}}e^{-\lvert F_{n} \rvert\beta+f_{F_n}(x_i)}\\
	&= e^{\lvert F_{n} \rvert \delta}\sum_{{i\in I}}e^{-\lvert F_{n} \rvert(\beta+\delta)+f_{F_n}(x_i)}\\ &\geq e^{\lvert F_{n} \rvert \delta}\sum_{{i\in I}}\mu(B_{F_n}(x_i,\epsilon))\\
	&\geq e^{\lvert F_{n} \rvert \delta}\mu(E_n)\\ &\geq \frac{e^{\lvert F_{n} \rvert \delta}}{n(n+1)}\mu(E).
	\end{align*}	
	Since $\lim_{n \to \infty}\frac{|F_n|}{\log n}=\infty $ and $\mu(E)>0,$ we have
	$$M^{P} (s,\frac{\epsilon}{10},E,\{F_n\}_{n=1}^{\infty},f)=\infty.$$
\end{proof}
\begin{prop}\label{prop 2.13}
	Let $\mu \in \mathcal{M} (X)$, $Z \subset X,$ $f \in C(X,\mathbb{R})$ and $\{F_n\}_{n=1}^{\infty}\subset\mathbb F(G)$ with $\lvert F_n\rvert \to \infty.$
	We have $$\overline P_\mu(\{F_n\}_{n=1}^{\infty},f) \leq P_\mu^{KP}(\{F_n\}_{n=1}^{\infty},f).$$
\end{prop}
\begin{proof}
For any $s < \overline P_\mu(\{F_n\}_{n=1}^{\infty},f),$ we can find a Borel set $A \subset X$ with $\mu(A)>0$ such that, for any $x \in A,$
$$\lim_{\epsilon \to 0}\limsup_{n \to \infty}\frac{-\log\mu(B_{F_n}(x,\epsilon))+f_{F_n}(x)}{\lvert F_{n} \rvert} >s.$$
Given $\delta \in (0,\mu(A))$ and $\epsilon>0.$ We shall show that $P_\mu^{KP}(\frac{\epsilon}{10},\delta,\{F_n\}_{n=1}^{\infty},f) > s,$ which implies that $P_\mu^{KP}(\{F_n\}_{n=1}^{\infty},f) \geq s.$
It suffices to show that
$$M_\mu^\mathcal{P}(s,\frac{\epsilon}{10},\delta,\{F_n\}_{n=1}^{\infty},f)=\infty.$$
Let $\{Z_i\}_{{i \in I}}$ be a finite or countable family with $\mu(\cup_{i \in I} Z_i)>1-\delta.$
Since
$$A=(A\cap (\cup_{i \in I}Z_i))\cup(A\verb|\|\cup_{i \in I}Z_i),$$
it follows that $\mu((A\cap (\cup_{i \in I}Z_i))) \geq \mu(A)-\delta>0.$ Thus there exists $i$ such that $\mu(A\cap Z_i)>0.$
Due to Proposition \ref{prop 2.12}, we have
$$M^{P}(s,\frac{\epsilon}{10},Z_i,\{F_n\}_{n=1}^{\infty},f) \geq M^{P}(s,\frac{\epsilon}{10},A\cap Z_i,\{F_n\}_{n=1}^{\infty},f)=\infty.$$
Thus 
$$M_\mu^\mathcal{P}(s,\frac{\epsilon}{10},\delta,\{F_n\}_{n=1}^{\infty},f)=\infty.$$


\end{proof}
\section{Proof of Theorem \ref{thm 1.1}}
In this section, we will prove the variational principle, which is divided into two parts: upper bound and lower bound.
\subsection{Lower bound}
Using Proposition \ref{prop 2.10} and Proposition \ref{prop 2.13}, this shows that
\begin{align*}
	&\sup\{\overline P_\mu(\{F_n\}_{n=1}^{\infty},f):\mu \in \mathcal{M} (X),\mu(Z)=1\}\\
	 \leq& \sup\{P_\mu^{KP}(\{F_n\}_{n=1}^{\infty},f):\mu \in \mathcal{M} (X),\mu(Z)=1\}\\
=& \sup\{P_\mu^P(\{F_n\}_{n=1}^{\infty},f):\mu \in \mathcal{M} (X),\mu(Z)=1\}\\
\leq &P^P(Z,\{F_n\}_{n=1}^{\infty},f).
\end{align*}
\subsection{Upper bound}
Now we just need to prove that for any non-empty analytic subset $Z \subset X,$ $f \in C(X,\mathbb{R})$ and $s \in ({\left\lVert f\right\rVert},P^P(Z,\{F_n\}_{n=1}^{\infty},f))$ there exists  $\mu$ satisfying $\mu(Z)=1$ and $\overline P_\mu(\{F_n\}_{n=1}^{\infty},f) \geq s.$

\begin{lem}\label{lem 3.1}
Let $Z \subset X,$ $\epsilon >0$ and $s >{\left\lVert f\right\rVert}. $ If $M^P(s,\epsilon,Z,\{F_n\}_{n=1}^{\infty},f)=\infty,$ then for any given finite interval
$(a,b) \subset [0,+\infty)$ and $N \in \mathbb{N}, $ there exists a finite disjoint collection $\{\overline B_{F_{n_i}}(x_i,\epsilon)\}$ such that $x_i \in Z,$ $n_i \geq N,$
and $\sum_{i}e^{-s\lvert F_{n_i} \rvert + f_{F_{n_i}}(x_i)} \in (a,b).$
\end{lem}
\begin{proof}
	Take $N_1>N$ large enough such that $e^{|F_{N_1}|(\left\lVert f\right\rVert-s)}<b-a.$ Since $M^P(s,\epsilon,Z,\{F_n\}_{n=1}^{\infty},f)=\infty,$ it follows that $M^P(N_1,s,\epsilon,Z,\{F_n\}_{n=1}^{\infty},f)=\infty.$
	There hence exists a finite disjoint collection $\{\overline B_{F_{n_i}}(x_i,\epsilon)\}$ such that $x_i \in Z,$ $n_i \geq N_1$ and $\sum_{i}e^{-|F_{n_i}|s+f_{{F_{n_i}}}(x_i)}>b.$
	Since $e^{-|F_{n_i}|s+f_{{F_{n_i}}}(x_i)} \leq e^{|F_{n_i}|(\left\lVert f\right\rVert-s)} \leq b-a,$ we can discard elements in this collection one by one until we have $\sum_{i}e^{-s\lvert F_{n_i} \rvert + f_{F_{n_i}}(x_i)} \in (a,b).$
\end{proof}
We now turn to show the Upper bound. We employ the approach used by Feng and Huang in \cite{re10}. For any $s \in ({\left\lVert f\right\rVert},P^P(Z,\{F_n\}_{n=1}^{\infty},f)),$ we take
 $\epsilon$ small enough such that $s < P^P(\epsilon,Z,\{F_n\}_{n=1}^{\infty},f).$ Fix $t \in (s,P^P(\epsilon,Z,\{F_n\}_{n=1}^{\infty},f)).$ 
Since $Z$ is analytic, there exists a continuous surjective map $\phi :\mathcal{N} \to Z.$ Let  $\Gamma _{n_1,n_2,...,n_p}=\{(m_1,m_2,...)\in\mathcal{N}:m_1 \leq n_1,m_2 \leq n_2,...,m_p \leq n_p \} $
and let $Z_{n_1,...,n_p}=\phi(\Gamma _{n_1,n_2,...,n_p}) .$
The construction is divided into the following  three steps:\\
\emph{Step 1.} Construct $K_1,$ $\mu_1,$ $n_1,$ $\gamma_1,$ and  $m_1(·) .$\\
Note that $M^\mathcal{P} (t,\epsilon,Z,\{F_n\}_{n=1}^{\infty},f)=\infty.$ Let
$$H=\bigcup\{G \subset X:G~is~open,M^\mathcal{P} (t,\epsilon,Z\cap G,\{F_n\}_{n=1}^{\infty},f)=0\}.$$
Then $M^\mathcal{P} (t,\epsilon,Z\cap H,\{F_n\}_{n=1}^{\infty},f)= 0$ by the separability of $X$. Let 
$$Z'=Z\verb|\|H=Z\cap(X\verb|\|H).$$
For any open set $G \subset X,$ either $Z'\cap G=\emptyset$ or $M^\mathcal{P} (t,\epsilon,Z'\cap G,\{F_n\}_{n=1}^{\infty},f)>0.$
Indeed, suppose that $M^\mathcal{P} (t,\epsilon,Z'\cap G,\{F_n\}_{n=1}^{\infty},f)=0.$
Since $Z=Z'\cup (Z\cap H),$
\begin{align*}
	M^\mathcal{P} (t,\epsilon,Z\cap G,\{F_n\}_{n=1}^{\infty},f) &\leq M^\mathcal{P} (t,\epsilon,Z'\cap G,\{F_n\}_{n=1}^{\infty},f)\\ &+M^\mathcal{P} (t,\epsilon,Z\cap H,\{F_n\}_{n=1}^{\infty},f)\\&=0.
\end{align*}
Thus $G \subset H,$ which implies that $Z'\cap G=\emptyset.$
Since
$$M^\mathcal{P} (t,\epsilon,Z,\{F_n\}_{n=1}^{\infty},f) \leq M^\mathcal{P} (t,\epsilon,Z\cap H,\{F_n\}_{n=1}^{\infty},f)+M^\mathcal{P} (t,\epsilon,Z',\{F_n\}_{n=1}^{\infty},f) $$
and $M^\mathcal{P} (t,\epsilon,Z\cap H,\{F_n\}_{n=1}^{\infty},f)=0,$ we have
$$M^\mathcal{P} (t,\epsilon,Z',\{F_n\}_{n=1}^{\infty},f)=M^\mathcal{P} (t,\epsilon,Z,\{F_n\}_{n=1}^{\infty},f)=\infty.$$
Using Lemma \ref{lem 3.1}, we can find a finite set $K_1 \subset Z',$ an integer-valued function $m_1(x)$ on $K_1$ such that the collection $ \{\overline B_{F_{m_1(x)}}(x,\epsilon)\}_{x \in K_1}$ is disjoint and
$$\sum_{x \in K_1}  e^{-s\lvert F_{m_1(x)} \rvert +f_{F_{m_1(x)}}(x)} \in (1,2).$$
Define
$$\mu_1=\sum_{x \in K_1}e^{-s\lvert F_{m_1(x)} \rvert +f_{F_{m_1(x)}}(x)}\delta_x.$$
Take a small $\gamma_1$ such that for any function $z:K_1 \to X$ with \\$\max_{x \in K_1}d(x,z(x))\leq  \gamma_1,$
we have for each $x \in K_1,$
$(\overline B(z(x),\gamma_1)\cup\overline B_{F_{m_1(x)}}(z(x),\epsilon))\cap (\bigcup_{y \in K_1\verb|\|\{x\}}\overline B(z(y),\gamma_1)\cup \overline B_{F_{m_1(y)}}(z(y),\epsilon))=\emptyset.$
It follows from $K_1 \subset Z'$ that for any $x \in K_1,$
\begin{align*}
	M^\mathcal{P} (t,\epsilon,Z \cap B(x,\frac{\gamma_1}{4}),\{F_n\}_{n=1}^{\infty},f) &\geq M^\mathcal{P} (t,\epsilon,Z' \cap B(x,\frac{\gamma_1}{4}),\{F_n\}_{n=1}^{\infty},f)\\&>0.
\end{align*}
Therefore we can pick a sufficiently large $n_1 \in \mathbb{N}$ so that $K_1 \subset Z_{n_1}$ and $M^\mathcal{P} (t,\epsilon,Z_{n_1} \cap B(x,\frac{\gamma_1}{4}),\{F_n\},f)>0$ for each $x \in K_1.$\\
\emph{Step 2.} Construct $K_2,$ $\mu_2,$ $n_2,$ $\gamma_2,$ and  $m_2(·) .$\\
The family of balls $\{\overline B(x,\gamma_1)\}_{x \in K_1}$ are pairwise disjoint. For each $x \in K_1,$ since $M^\mathcal{P} (t,\epsilon,Z_{n_1} \cap B(x,\frac{\gamma_1}{4}),\{F_n\},f)>0,$
we can construct a finite set  as in  Step 1  
$$E_2(x) \subset Z_{n_1} \cap B(x,\frac{\gamma_1}{4})$$
and an integer-valued function
$$m_2:E_2(x) \to \mathbb{N}\cap[\max\{m_1(y):y \in K_1\},\infty]$$
such that
\begin{flushleft}
(2-a)$M^\mathcal{P} (t,\epsilon,Z_{n_1} \cap G,\{F_n\},f)>0 ,$ for any open set $G$ with $G\cap E_2(x) \neq \emptyset;$\\
(2-b)the elements in $\{\overline B_{F_{m_2(y)}}(y,\epsilon)\}_{y \in E_2(x)}$ are disjoint, and
\end{flushleft}
$$\mu_1(\{x\})<\sum_{y \in E_2(x)}e^{-s\lvert F_{m_2(y)} \rvert +f_{F_{m_2(y)}}(y)}<(1+2^{-2})\mu_1(\{x\}).$$
To see it, we fix $x \in K_1.$
Denote $F=Z_{n_1} \cap B(x,\frac{\gamma_1}{4}).$ 
Let
$$H_x=\bigcup\{G \subset X:G~is~open,M^\mathcal{P} (t,\epsilon,F\cap G,\{F_n\}_{n=1}^{\infty},f)=0\}.$$
Set 
$$F'=F\verb|\|H_x.$$
Then as in {Step 1,} we can show that
$$M^\mathcal{P} (t,\epsilon,F',\{F_n\}_{n=1}^{\infty},f)=M^\mathcal{P} (t,\epsilon,F,\{F_n\}_{n=1}^{\infty},f)>0$$
and
$$M^\mathcal{P} (t,\epsilon,F'\cap G,\{F_n\}_{n=1}^{\infty},f)>0,$$
for any open set $G$ with $G \cap F'\neq \emptyset.$ Since $s < t,$
$$M^\mathcal{P} (s,\epsilon,F',\{F_n\}_{n=1}^{\infty},f)=\infty.$$
Using Lemma \ref{lem 3.1} again, we can find a finite set $ E_2(x) \subset F',$ an integer-valued function $m_2(x)$ on $E_2(x)$ so that (2-b) holds.
Observe that if $G \cap E_2(x)\neq \emptyset$ and $G$ is open, then $G \cap F'\neq \emptyset.$ Hence
$$M^\mathcal{P} (t,\epsilon,Z_{n_1}\cap G,\{F_n\}_{n=1}^{\infty},f) \geq M^\mathcal{P} (t,\epsilon,F'\cap G,\{F_n\}_{n=1}^{\infty},f)>0.$$
Then (2-a) holds.
Since the family $\{\overline{B}(x,\gamma_1)\}_{x \in K_1}$ is disjoint, $E_2(x) \cap E_2(x')=\emptyset$ for different $x,x' \in K_1.$
Define $K_2=\bigcup_{x \in K_1}E_2(x)$
and
$$\mu_2=\sum_{y \in K_2}e^{-s\lvert F_{m_2(y)} \rvert +f_{F_{m_2(y)}}(y)}\delta_y.$$
The elements in  $\{\overline B_{F_{m_2(y)}}(y,\epsilon)\}_{y \in K_2} $ are disjoint. Hence we can take $\gamma_2 \in (0,\frac{\gamma_1}{4})$ such that
for any function $z:K_2 \to X$ satisfying \\$\max_{x \in K_2}d(x,z(x))<\gamma_2,$
we have
$(\overline B(z(x),\gamma_2)\cup \overline B_{F_{m_2(x)}}(z(x),\epsilon))\cap (\bigcup_{y \in K_2\verb|\|{x}}\overline B(z(y),\gamma_2)\cup \overline B_{F_{m_2(y)}}(z(y),\epsilon))=\emptyset.$
For each $x \in K_2.$
Choose a sufficiently large $n_2 \in \mathbb{N}$ so that $K_2 \subset Z_{n_1,n_2}$ and 
$$M^\mathcal{P} (t,\epsilon,Z_{n_1,n_2} \cap B(x,\frac{\gamma_2}{4}),\{F_n\}_{n=1}^{\infty},f)>0,$$
for each $x \in K_2.$\\
\emph{Step 3.} Assume that $K_i,$ $\mu_i,$ $n_i,$ $\gamma_i,$ and  $m_i(·)$ have been constructed for $i=1,...,p.$
In particular, suppose that for any function $z:K_p \to X$ with $\max_{x \in K_p}d(x,z(x))<\gamma_p$,
we have
$(\overline B(z(x),\gamma_p)\cup \overline B_{F_{m_p(x)}}(z(x),\epsilon))\cap (\bigcup_{y \in K_p\verb|\|\{x\}}\overline B(z(y),\gamma_p)\cup \overline B_{F_{m_p(y)}}(z(y),\epsilon))=\emptyset,$
for each $x \in K_p.$
Then $K_p \subset Z_{n_1...n_p}$ and 
$$M^\mathcal{P} (t,\epsilon,Z_{n_1...n_p} \cap B(x,\frac{\gamma_p}{4}),\{F_n\}_{n=1}^{\infty},f)>0,$$
for each $x \in K_p.$
The family of balls $\{\overline B(x,\gamma_p)\}_{x \in K_p} $ are pairwise disjoint. For each $x \in K_p,$ since $M^\mathcal{P} (t,\epsilon,Z_{n_1...n_p} \cap B(x,\frac{\gamma_p}{4}),\{F_n\}_{n=1}^{\infty},f)>0,$
we can construct as in Step 2 a finite set 
$$E_{p+1}(x) \subset Z_{n_1...n_p} \cap B(x,\frac{\gamma_p}{4})$$
and an integer-valued function
$$m_{p+1}:E_{p+1}(x) \to \mathbb{N}\cap[\max\{m_p(y):y \in K_p\},\infty]$$
such that
\begin{flushleft}
(3-a)$M^\mathcal{P} (t,\epsilon,Z_{n_1...n_p} \cap G,\{F_n\},f)>0,$ for any open set $G$ with $G\cap E_{p+1}(x) \neq \emptyset;$\\
(3-b)the elements in $\{\overline B_{F_{m_{p+1}(y)}}(y,\epsilon)\}_{y \in E_{p+1}(x)}$ are disjoint, and
\end{flushleft}
$$\mu_p(\{x\})<\sum_{y \in E_{p+1}(x)}e^{-s\lvert F_{m_{p+1}(y)} \rvert +f_{F_{m_{p+1}(y)}}(y)}<(1+2^{-p-1})\mu_p(\{x\}).$$
Clearly  $E_{p+1}(x) \cap E_{p+1}(y) = \emptyset$
for different $x,y \in K_p.$
Define $K_{p+1}=\bigcup_{x \in K_p}E_{p+1}(x)$
and
$$\mu_{p+1}=\sum_{y \in K_{p+1}}e^{-s\lvert F_{m_{p+1}(y)} \rvert +f_{F_{m_{p+1}(y)}}(y)}\delta_y.$$
The elements in  $\{\overline B_{F_{m_{p+1}(y)}}(y,\epsilon)\}_{y \in K_{p+1}} $ are disjoint. Hence we can take $\gamma_{p+1} \in (0,\frac{\gamma_p}{4})$ such that
for any function $z:K_{p+1} \to X$ satisfying $\max_{x \in K_{p+1}}d(x,z(x))<\gamma_{p+1},$
we have
$(\overline B(z(x),\gamma_{p+1})\cup \overline B_{F_{m_{p+1}(x)}}(z(x),\epsilon)) \\ \cap (\bigcup_{y \in K_{p+1}\verb|\|\{x\}}\overline B(z(y),\gamma_{p+1})\cup \overline B_{F_{m_{p+1}(y)}}(z(y),\epsilon))=\emptyset,$
for each $x \in K_{p+1}.$
Choose a sufficiently large $n_{p+1} \in \mathbb{N}$ so that $K_{p+1} \subset Z_{n_1...n_{p+1}}$ and 
$$M^\mathcal{P} (t,\epsilon,Z_{n_1...n_{p+1}} \cap B(x,\frac{\gamma_{p+1}}{4}),\{F_n\}_{n=1}^{\infty},f)>0,$$
for each $x \in K_{p+1}.$
As in above steps, we can construct by introduction $\{K_i\},\{\mu_i\},n_i,\gamma_i,and~m_i(·) .$
We summarize some of their basic properties as follows:\\
(a)For each $i$, the family $\mathcal{F}_i:=\{\overline{B}(x,\gamma_i):x \in K_i\}$ is disjoint. For every $B \in \mathcal{F}_{i+1},$ there exists $x \in K_i$ such that $B \subset {B}(x,\frac{\gamma_i}{2}).$ \\
(b)For each $x \in K_i$ and $z\in {B}(x,\gamma_i),$ we have
$$\overline B_{F_{m_{i}(x)}}(z,\epsilon) \cap \bigcup_{y \in K_i\verb|\|\{x\}}  B(y,\gamma_i)=\emptyset$$
and
\begin{align*}
	\mu_i( B(x,\gamma_i))&=e^{-s\lvert F_{m_i(x)} \rvert+f_{F_{m_i(x)}(x)}}\\ &<\sum_{y \in E_{i+1}(x)}e^{-s\lvert F_{m_{i+1}(y)} \rvert +f_{F_{m_{i+1}(y)}}(y)}\\
&<(1+2^{-i-1})\mu_i( B(x,\gamma_i)),
\end{align*}
where $E_{i+1}(x)=B(x,\gamma_i) \cap K_{i+1}.$
Furthermore, for $F_i\in \mathcal{F}_{i}$,
$$\mu_i(F_i) \leq \mu_{i+1}(F_i)=\sum_{F\in \mathcal{F}_{i+1}:F \subset F_i }\mu_{i+1}(F)$$
$$\leq \sum_{F\in \mathcal{F}_{i+1}:F \subset F_i }(1+2^{-i-1})\mu_i(F)$$
$$=(1+2^{-i-1})\sum_{F\in \mathcal{F}_{i+1}:F \subset F_i }\mu_i(F)$$
$$\leq (1+2^{-i-1})\mu_i(F_i).$$ 
Using the above inequalities repeatedly, we have for any $j>i,$ $F_i \in \mathcal{F}_i,$
\begin{align}\label{3.1}
\mu_i(F_i) \leq \mu_j(F_i) \leq \prod _{n=i+1}^j (1+2^{-n})\mu_i(F_i) \leq C\mu_i(F_i),
\end{align}
where $C:=\prod _{n=1}^\infty(1+2^{-n})<\infty.$
Let $\widetilde\mu$ be a limit point of $\{\mu_i\}$ in the weak-star topology, let
$$K=\bigcap_{n=1}^{\infty}\overline{\bigcup_{i \geq n}K_i}.$$
Then $\widetilde\mu$ is supported on $K,$ $K \subset \bigcap_{p=1}^{+\infty}\overline{Z_{n_1,...,n_p}}.$
By the continuity of $\phi,$ applying the Cantor's diagonal argument, we can show that $\bigcap_{p=1}^{+\infty}\overline{Z_{n_1,...,n_p}} =\bigcap_{p=1}^{+\infty}{Z_{n_1,...,n_p}}.$
Hence $K$ is a compact subset of $Z$. For any $x \in K_i,$ by \ref{3.1} we have 
\begin{align*}
	e^{-s\lvert F_{m_i(x)} \rvert+f_{F_{m_i(x)}(x)}}&=\mu_i( B(x,\gamma_i))\\&\leq \widetilde\mu( B(x,\gamma_i))\\
&\leq C\mu_i( B(x,\gamma_i))\\&=Ce^{-s\lvert F_{m_i(x)} \rvert+f_{F_{m_i(x)}(x)}}.
\end{align*}
In particular,
\begin{align*}
	1\leq \sum_{x \in K_1}\mu_1(B(x,\gamma_1)) 
&\leq\widetilde\mu(K)\\ &\leq \sum_{x \in K_1}C\mu_1(B(x,\gamma_1))\\&\leq 2C.
\end{align*}
For every $x \in K_i$ and $z \in \overline{B}(x,\gamma_i),$
$$\widetilde\mu(\overline B_{F_{m_i(x)}}(z,\epsilon))\leq \widetilde \mu(\overline B(x,\frac{\gamma_i}{2})) \leq Ce^{-s\lvert F_{m_i(x)} \rvert+f_{F_{m_i(x)}(x)}}.$$
For each $z \in K$ and $i \in \mathbb{N},$ $z \in \overline B(x,\frac{\gamma_i}{2})$ for some $x \in K_i$. Thus
$$\widetilde\mu(\overline B_{F_{m_i(x)}}(z,\epsilon)) \leq Ce^{-s\lvert F_{m_i(x)} \rvert+f_{F_{m_i(x)}(x)}}.$$
Let $\mu=\widetilde\mu \verb|\| \widetilde\mu(K).$ Then $\mu \in \mathcal{M} (X),$ $\mu(K)=1,$ and for every $z \in K,$ 
there exists a sequence $\{k_i\}$ with $k_i \to \infty$ such that
$$\mu(B_{F_{k_i}}(z,\epsilon)) \leq \frac{Ce^{-s\lvert F_{k_i} \rvert+f_{F_{k_i}(z)}}}{\widetilde\mu(K)}.$$
This implies that $\overline P_\mu(\{F_n\}_{n=1}^{\infty},f) \geq s.$

\section{Proof of Theorem \ref{thm 1.2}}
\subsection{Proof of the first inequality}
We need to prove $$P^P(\pi(E),\{F_n\}_{n=1}^{\infty},f) \leq P^P(E,\{F_n\}_{n=1}^{\infty},f\circ \pi).$$
Let $G,\{F_n\}_{n=1}^{\infty}$ and $\pi:(X,G) \to (Y,G)$ be as in Theorem \ref{thm 1.2} and $E\subset X$ be a subset.
Let $d$ and $\rho $ be the compatible metrics on $X$ and $Y,$ respectively.
For any $\epsilon >0 ,$ there exists $\delta >0$ such that for any $x_1,x_2 \in X$ with $d(x_1,x_2 ) \leq \delta,$
one has $\rho(\pi(x_1),\pi(x_2)) \leq \epsilon.$
Now let $\{y_i\}_{i=1}^k \subset \pi(E)$  be any $(F_n,\epsilon)$ separated set of $\pi(E) $ and  for each $i$ choose a point $x_i \in \pi^{-1}(y_i)\cap E.$
Hence
\begin{align*}
	\sum_{i=1}^k e^{f_{F_n}(y_i)} &= \sum_{i=1}^k e^{f_{F_n} \circ \pi (x_i)} \\
	&\leq P(E,\delta,F_n,f\circ \pi).
\end{align*}
This implies 
$$P(\pi(E),\epsilon,F_n,f) \leq P(E,\delta,F_n,f\circ \pi).$$
Furthermore,
$$P^{UC}(\pi(E),\epsilon,\{F_n\}_{n=1}^{\infty},f ) \leq P^{UC}(E,\delta,\{F_n\}_{n=1}^{\infty},f\circ \pi).$$
By Proposition \ref{prop 2.5}, for any $\eta >0,$ there exists a cover $\bigcup_{i=1}^\infty E_i=E$ such that
$$P^{P}(E,\frac{\delta}{3},\{F_n\}_{n=1}^{\infty},f \circ \pi) + \eta \geq \sup_{i \geq 1}P^{UC}(E_i,\delta,\{F_n\}_{n=1}^{\infty},f\circ \pi).$$
Then we have
\begin{align*}
	P^{P}(\pi(E),\epsilon,\{F_n\}_{n=1}^{\infty},f ) &\leq \sup_{i \geq 1}P^{P}(\pi(E_i),\epsilon,\{F_n\}_{n=1}^{\infty},f)\\
&\leq \sup_{i \geq 1}P^{UC}(\pi(E_i),\epsilon,\{F_n\}_{n=1}^{\infty},f )\\
&\leq \sup_{i \geq 1}P^{UC}(E_i,\delta,\{F_n\}_{n=1}^{\infty},f \circ\pi)\\
&\leq P^{P}(E,\frac{\delta}{3},\{F_n\}_{n=1}^{\infty},f \circ \pi) + \eta,
\end{align*}
which implies that 
$$P^P(\pi(E),\{F_n\}_{n=1}^{\infty},f) \leq P^P(E,\{F_n\}_{n=1}^{\infty},f\circ \pi).$$
\subsection{Proof of the second inequality}

\begin{lem}\label{lem 4.1}(A claim in \cite{re1})
	Fix $\tau>0$, for any $\epsilon>0$. When $\eta$ is small enough and $N$ is large enough, for all $y \in Y$ and $n >N$, there exists $l(y)>0$ and $v_1(y),...,v_{l(y)}(y) \in X$ such that
	$$\pi^{-1}(B_{F_n}(y,\eta,\rho))    \subset    \bigcup_{i=1}^{l(y)}B_{F_n}(v_i(y),4\epsilon,d),$$
	$$l(y) \leq \exp((a+2\tau)\lvert F_n \rvert)$$
	and
	$$ \pi B_{F_n}(v_i(y),4\epsilon,d) \cap B_{F_n}(y,\eta,\rho)\neq \emptyset,1 \leq i \leq l(y),$$
	where $a=\sup_{y \in Y}h_{top}^{UC}(\pi^{-1}(y),\{F_n\}_{n=1}^{\infty})<\infty.$
\end{lem}

\begin{thm}\label{thm 4.2}
	Let G be a countable infinite discrete amenable group, $\pi:(X,G)\to (Y,G)$ be a factor map, $\{F_n\}_{n=1}^{\infty}$ be any tempered F$\o$lner sequence   satisfying $\lim_{n \to \infty}\frac{|F_n|}{\log n}=\infty, E \subset X $  and $f \in C(X,\mathbb{R})$.
	Then 
	$$ P^{UC}(E,\{F_n\}_{n=1}^{\infty},f\circ \pi) \leq P^{UC}(\pi(E),\{F_n\}_{n=1}^{\infty},f)+\sup_{y \in Y}h_{top}^{UC}(\pi^{-1}(y),\{F_n\}_{n=1}^{\infty}),$$
	where $ h_{top}^{UC}(Z,\{F_n\}_{n=1}^{\infty})=P^{UC}(Z,\{F_n\}_{n=1}^{\infty},0).$
\end{thm}
\begin{proof}
	We assume that
	$$a:=\sup_{y \in Y}h_{top}^{UC}(\pi^{-1}(y),\{F_n\}_{n=1}^{\infty})<\infty.$$
	Fix $\tau >0$, $\epsilon>0$ and $\eta>0 $ small enough. Let $N$ be as in Lemma \ref{lem 4.1}. For any $n >N,$
    we will apply Lemma \ref{lem 4.1} in the following step.
	For any subset $E$ of $X$, let $H$ be an $(F_n,\eta)$-spanning set of $\pi(E)$ with minimal cardinality. Then by the claim, the set
	$R=\{v_i(y):1 \leq i \leq l(y) ,y \in H\}$ forms an $(F_n,4\epsilon)$-spanning set of $E$. Since
	$$E \subset \pi^{-1}\pi(E) \subset \bigcup_{y \in H}\pi^{-1}B_{F_n}(y,\eta,\rho) \subset  \bigcup_{y \in H}\bigcup_{i=1}^{l(y)}B_{F_n}(v_i(y),4\epsilon,d).$$
    Assume that$~var(f,\eta)=\sup_{x,y \in X}\{|f(x)-f(y)|:d(x,y) \leq \eta\}$ and $var(f\circ \pi,4\epsilon)=\sup_{x,y \in Y}\{|f \circ \pi(x)-f \circ \pi(y)|:\rho(x,y) \leq 4\epsilon\}.$
    Thus
	$$ Q(E,4\epsilon,F_n,f \circ \pi) \leq \sum_{y \in H} \sum_{i=1}^{l(y)} \exp(f_{F_n} \circ \pi(v_i(y)))$$
	$$\leq \sum_{y \in H} \sum_{i=1}^{l(y)} \exp(f_{F_n}(y) +\lvert F_n \rvert var(f,\eta)+\lvert F_n \rvert  var(f\circ \pi,4\epsilon))$$
	$$\leq \exp((a+2\tau)\lvert F_n \rvert)\sum_{y \in H} \exp(f_{F_n}(y) +\lvert F_n \rvert var(f,\eta)+\lvert F_n \rvert  var(f\circ \pi,4\epsilon)).$$
	Thus
	$$\limsup_{n \to \infty}  \frac{1}{\lvert F_n \rvert} \log Q(E,4\epsilon,F_n,f \circ \pi) $$
	$$\leq a+2\tau+var(f\circ \pi,4\epsilon)+var(f,\eta)+\limsup_{n \to \infty}  \frac{1}{\lvert F_n \rvert} \log\sum_{y \in H} \exp(f_{F_n}(y)).$$
	$P^{UC'}(E,4\epsilon,\{F_n\}_{n=1}^{\infty},f\circ \pi) \leq P^{UC'}(\pi(E),\eta,\{F_n\}_{n=1}^{\infty},f)+a+2\tau+var(f\circ \pi,4\epsilon)+var(f,\eta).$
	Letting $\epsilon \to 0$, $\eta \to 0$ and $\tau \to 0,$ we have
	$$P^{UC}(E,\{F_n\}_{n=1}^{\infty},f\circ \pi) \leq P^{UC}(\pi(E),\{F_n\}_{n=1}^{\infty},f)+\sup_{y \in Y}h_{top}^{UC}(\pi^{-1}(y),\{F_n\}_{n=1}^{\infty}).$$

\end{proof}
Next we prove the the second inequality.
By Theorem \ref{thm 4.2}, we have 
\begin{align*}
&P^{UC}(E,8\epsilon,\{F_n\}_{n=1}^{\infty},f\circ \pi)-var(f\circ \pi,4\epsilon)\\
&\leq P^{UC'}(E,4\epsilon,\{F_n\}_{n=1}^{\infty},f\circ \pi) \\
&\leq P^{UC'}(\pi(E),\eta,\{F_n\}_{n=1}^{\infty},f)+a+2\tau+var(f\circ \pi,4\epsilon)+var(f,\eta)\\
&\leq P^{UC}(\pi(E),\eta,\{F_n\}_{n=1}^{\infty},f)+a+2\tau+var(f\circ \pi,4\epsilon)+var(f,\eta).	
\end{align*}

By Proposition \ref{prop 2.5}(2), for all $\delta>0,$ there exists a cover $\bigcup_{i=1}^{\infty}V_i=\pi(E)$ such that
$$P^P(\pi(E),\frac{\eta}{3},\{F_n\}_{n=1}^{\infty},f) +\delta \geq \sup_{i \geq 1}P^{UC}(V_i,\eta,\{F_n\}_{n=1}^{\infty},f).$$
Therefore,
\begin{align*}
&P^{P}(E,8\epsilon,\{F_n\}_{n=1}^{\infty},f\circ \pi) \leq \sup_{i \geq 1}P^{P}(\pi^{-1}(V_i),8\epsilon,\{F_n\}_{n=1}^{\infty},f\circ \pi) \\
	&\leq \sup_{i \geq 1}P^{UC}(\pi^{-1}(V_i),8\epsilon,\{F_n\}_{n=1}^{\infty},f\circ \pi)\\
	&\leq \sup_{i \geq 1}P^{UC}(V_i,\eta,\{F_n\}_{n=1}^{\infty},f)+a+2\tau+2var(f\circ \pi,4\epsilon)+var(f,\eta)\\
	&\leq P^P(\pi(E),\frac{\eta}{3},\{F_n\}_{n=1}^{\infty},f) +\delta +a+2\tau+2var(f\circ \pi,4\epsilon)+var(f,\eta).
\end{align*}
Hence,
$$P^{P}(E,\{F_n\}_{n=1}^{\infty},f\circ \pi) \leq P^P(\pi(E),\{F_n\}_{n=1}^{\infty},f)+a+2\tau.$$
Since $\tau$ is arbitrary, we finally obtain
$$P^P(E,\{F_n\}_{n=1}^{\infty},f\circ \pi) \leq P^P(\pi(E),\{F_n\}_{n=1}^{\infty},f)+\sup_{y \in Y}h_{top}^{UC}(\pi^{-1}(y),\{F_n\}_{n=1}^{\infty}).$$
\section{Proof of Theorem \ref{thm 1.3}}	

In the first subsection, we will apply a lemma in \cite{re3} to show that the upper bound which does not require the properties of ergodic and almost $specification ~property$.
In the second and third subsections, we will prove the lower bound under the {\it ergodic assumption or almost  specification  property}, respectively.
Zhang proved the case of almost  specification  property  for Pesin-Pitskel pressure in \cite{re3}. We only focus on proving the case of ergodic.
\subsection{Proof of upper bound}
\begin{lem}\label{lem 5.1} \cite[Lemma 4.2]{re3}
	Let $\{F_n\}_{n=1}^{\infty}$ be a F$\o$lner sequence, $\mu \in \mathcal{M} (X,G)$, $C \subset \mathcal{M} (X)$ be a neighborhood of $\mu,$ $\phi \in C(X,\mathbb{R})$ and set 
	$$P(X_{F_n,C},\epsilon,F_n,\phi)=\sup_{E_{F_n,C}}\sum_{x \in E_{F_n,C}}e^{\phi_{F_n}(x)},$$
    where the supremum is taken over all $(F_n,\epsilon)$-separated sets $E_{F_n,C} \subset X_{F_n,C}$ and $X_{F_n,C}:=\{x \in X:\frac{1}{\lvert F_n \rvert}\sum_{g \in F_n }\delta_x \circ g^{-1} \in C\}.$
	Then  
	$$\lim_{\epsilon \to 0} \inf_{\{C:\mu \in C\}}\limsup_{n \to \infty}  \frac{1}{\lvert F_n \rvert} \log P(X_{F_n,C},\epsilon,F_n,\phi) \leq h_\mu(X,G)+\int \phi                                                                                                                                                                                                                                                                                                                                                                                                                                                                                                                                                                                                                                                                                                                                                                                                                                                              \,d\mu.$$
\end{lem}
In the following we are going to prove $P^P(G_\mu,\{F_n\},f) \leq h_\mu(X,G)+\int f \,d\mu$ assuming that the F$\o$lner sequence $\{F_n\}$ satisfies the 
growth condition(1.1).
For $\mu \in \mathcal{M}  (X,G) ,$ let $\{K_m\}_{m \in \mathbb{N}}$ be a decreasing sequence of closed convex neighborhoods of $\mu$ in $\mathcal{M}  (X)$
such that $\bigcap_{m \in \mathbb{N}}K_m=\{\mu\}.$
Let 
$$A_{n,m}=\{x\in X:\frac{1}{\lvert F_n \rvert}\sum_{g \in F_n }\delta_x \circ g^{-1} \in K_m \}for~m,n \in \mathbb{N},$$
and
$$R_{N,m}=\{x\in X:\forall n \geq N,\frac{1}{\lvert F_n \rvert}\sum_{g \in F_n }\delta_x \circ g^{-1} \in K_m \}for~m,N \in \mathbb{N}.$$
Then for any $m,N \geq 1,$
$$R_{N,m}=\bigcap_{n>N}A_{n,m}~and~G_\mu \subset \bigcup_{k>N}R_{k,m}.$$
For any $\epsilon >0$, we have
\begin{align*}
 &\limsup_{n \to \infty}  \frac{1}{\lvert F_n \rvert} \log P(R_{N,m},\epsilon,F_n,f) \\ 
 \leq &\limsup_{n \to \infty}  \frac{1}{\lvert F_n \rvert} \log P(A_{n,m},\epsilon,F_n,f), \text{ for any } m,N \geq 1. 
\end{align*}

By the Lemma \ref{lem 5.1},
$$\lim_{\epsilon \to 0}\lim_{m \to \infty}\limsup_{n \to \infty}  \frac{1}{\lvert F_n \rvert} \log P(A_{n,m},\epsilon,F_n,f) \leq h_\mu(X,G)+\int f \,d\mu.$$
Hence for any $\eta>0,$ there exists $\epsilon_1>0$ such that, for any $0<\epsilon<\epsilon_1,$ there exists $M=M(\epsilon)\in \mathbb{N}$ such that
$$\limsup_{n \to \infty}  \frac{1}{\lvert F_n \rvert} \log P(A_{n,m},\epsilon,F_n,f) \leq h_\mu(X,G)+\int f \,d\mu+\eta,$$
whenever $m \geq M.$ Especially,
$$\limsup_{n \to \infty}  \frac{1}{\lvert F_n \rvert} \log P(A_{n,M},\epsilon,F_n,f) \leq h_\mu(X,G)+\int f \,d\mu+\eta.$$
Hence for any $0<\epsilon<\epsilon_1,$ we have for any $N \in \mathbb{N},$
$$\limsup_{n \to \infty}  \frac{1}{\lvert F_n \rvert} \log P(R_{N,M},\epsilon,F_n,f) \leq h_\mu(X,G)+\int f \,d\mu+\eta.$$
Since for any $N' \in N,G_\mu \subset \bigcup_{N>N'}R_{N,M},$ we have
\begin{align*}
	P^P(G_\mu,\epsilon,\{F_n\}_{n=1}^{\infty},f) \leq &\sup_{N>N'}P^P(R_{N,M},\epsilon,\{F_n\}_{n=1}^{\infty},f)\\
	\leq &\sup_{N>N'}\limsup_{n \to \infty}  \frac{1}{\lvert F_n \rvert} \log P(R_{N,M},\epsilon,F_n,f)\\
	\leq &h_\mu(X,G)+\int f \,d\mu+\eta.
\end{align*}
Letting $\epsilon \to 0$ and $\eta \to 0,$ we finish  the proof.
\subsection{Proof of lower bound if $(X,G)$ satisfies almost specification property}
When $(X,G)$ satisfies almost specification ~property, by the Theorem 1.1 in \cite{re3}, we have
$$P^B(G_\mu,\{F_n\}_{n=1}^{\infty},f)=h_\mu(X,G)+\int f \,d\mu.$$
Since $P^P(G_\mu,\{F_n\}_{n=1}^{\infty},f) \geq P^B(G_\mu,\{F_n\}_{n=1}^{\infty},f),$
we obtain $$P^P(G_\mu,\{F_n\}_{n=1}^{\infty},f) \geq h_\mu(X,G)+\int f \,d\mu.$$

\subsection{Proof of lower bound if $\mu$ is ergodic and $\{F_n\}$ is tempered}
Let $(X,G),$ $\mu$ and $\{F_n\}_{n=1}^{\infty}$ be as in Theorem \ref{thm 1.3}. If $Y\subset X,$ $f \in C(X,\mathbb{R}),$ and $\mu(Y)=1$ we will show that
$$P^P(Y,\{F_n\}_{n=1}^{\infty},f) \geq h_\mu(X,G)+\int f \,d\mu.$$
Using \cite{re5} Theorem 2.1(Brin-Katok entropy formula: ergodic case), we have for $\mu$ almost everywhere $x \in X,$
$$\lim_{\epsilon \to 0}\limsup_{n \to \infty}\dfrac{-\log\mu(B_{F_n}(x,\epsilon))}{\lvert F_n \rvert}=h_\mu(X,G).$$
Then we have for $\mu$ almost everywhere $x \in X,$
$$\overline P_\mu(x,\{F_n\}_{n=1}^{\infty},f)=h_\mu(X,G)+\int f \,d\mu.$$
There exists $Y \subset X,$ $\mu(Y)=1$ satisfying $$\overline P_\mu(x,\{F_n\}_{n=1}^{\infty},f)=h_\mu(X,G)+\int f \,d\mu$$ for all $x \in Y.$
Using proposition \ref{prop 2.12}(2), we have $$P^P(G_\mu \cap Y,\{F_n\}_{n=1}^{\infty},f) \geq h_\mu(X,G)+\int f \,d\mu.$$  Hence
$$P^P(G_\mu,\{F_n\}_{n=1}^{\infty},f) \geq h_\mu(X,G)+\int f \,d\mu.$$

\section*{Acknowledgement} 
 The  second author was supported by the
National Natural Science Foundation of China (No.12071222).
The third authors was  supported by the
National Natural Science Foundation of China (No. 11971236), Qinglan Project of Jiangsu Province of China.  The work was also funded by the Priority Academic Program Development of Jiangsu Higher Education Institutions.  Besides, we would like to express our gratitude to Tianyuan Mathematical Center in Southwest China(No. 11826102), Sichuan University and Southwest Jiaotong University for their support and hospitality.

\section*{Data availability} 
No data was used for the research described in the article.
\section*{Conflict of interest} 
The author declares no conflict of interest.


\begin{thebibliography}{HD82}
\normalsize


\bibitem[1]{re21} R. L. Adler, A. G. Konheim, and M. H. McAndrew. Topological entropy. \emph{Trans. Amer. Math. Soc.} 114 (1965), 309–319.

\bibitem[2]{re17} R. Bowen. Entropy for group endomorphisms and homogeneous spaces. \emph{Trans. Amer. Math. Soc.} 153 (1971), 401–414.

\bibitem[3]{re8} R. Bowen. Topological entropy for noncompact sets. \emph{Trans. Amer. Math. Soc.} 184 (1973), 125–136.

\bibitem[4]{re13} M. Brin and A. Katok. On local entropy. \emph{Springer-Verlag, Berlin,} 1983, 30–38.

\bibitem[5]{re27} C. Fang, W. Huang, Y. Yi and P. Zhang. Dimensions of stable sets and scrambled sets in positive finite entropy systems. \emph{Ergodic Theory Dynam. Systems} 32 (2012), no. 2, 599–628.

\bibitem[6]{re22} E. I. Dinaburg. A correlation between topological entropy and metric entropy. \emph{Dokl. Akad. Nauk SSSR} 190 (1970), 19–22.

\bibitem[7]{re1} D. Dou, D. Zheng and X. Zhou. Packing topological entropy for amenable group actions. \emph{Ergodic Theory Dynam. Systems} 43 (2023), no. 2, 480–514.

\bibitem[8]{re10} D. J. Feng and W. Huang. Variational principles for topological entropies of subsets. \emph{J. Funct. Anal.} 263 (2012), no. 8, 2228–2254.

\bibitem[9]{re26} X. Huang, Z. Li and Y. Zhou. A variational principle of topological pressure on subsets for amenable group actions. \emph{Discrete Contin. Dyn. Syst.} 40 (2020), no. 5, 2687–2703.

\bibitem[10]{re14} D. Kerr and H. Li. Ergodic Theory: Independence and Dichotomies. \emph{Springer, Cham,} 2016, xxxiv+431 pp.

\bibitem[11]{re20} A. Kolomogorov. A new metric invariant of transient dynamical systems and automorphisms of lebesgue spaces. \emph{Dokl. Akad. Nauk SSSR} (N.S.)119(1958), 861–864.

\bibitem[12]{re18} E. Lindenstrauss. Pointwise theorems for amenable groups. \emph{Invent. Math.} 146 (2001), no. 2, 259–295.

\bibitem[13]{re28} Q. Li, E. Chen and X. Zhou.  Corrigendum to: 'A note of topological pressure for non-compact sets of a factor map.' \emph{Chaos Solitons Fractals} 53 (2013), 75–77.
                   
\bibitem[14]{re40} P. Mattila. \emph{Geometry of Sets and Measures in Euclidean Spaces.} Cambridge University Press, Cambridge, 1995.

\bibitem[15]{re31} P. Oprocha and G. Zhang. Dimensional entropy over sets and fibres. \emph{Nonlinearity} 24 (2011), no. 8, 2325–2346.

\bibitem[16]{re15} D. S. Ornstein and B. Weiss. Entropy and isomorphism theorems for actions of amenable groups. \emph{J. Analyse Math.} 48 (1987), 1–141.

\bibitem[17]{re12} Y. B. Pesin. Dimension theory in dynamical systems: contemporary views and applications. \emph{University of Chicago Press,} 2008.  

\bibitem[18]{re9} Y. B. Pesin and B. S. Pitskel. Topological pressure and the variational principle for noncompact sets. \emph{Funktsional. Anal. i Prilozhen.} 18 (1984), no. 4, 50–63, 96.

\bibitem[19]{re19} C. E. Pfister and W. G. Sullivan. On the topological entropy of saturated sets. \emph{Ergodic Theory Dynam. Systems} 27 (2007), no. 3, 929–956.

\bibitem[20]{re25} X. Tang, W. C. Cheng and Y. Zhao. Variational principle for topological pressures on subsets. \emph{J. Math. Anal. Appl.} 424 (2015), no. 2, 1272–1285.

\bibitem[21]{re16} P. Walters. An introduction to ergodic theory, \emph{Springer-Verlag, New York-Berlin,} 1982, ix+250 pp.

\bibitem[22]{re24} C. Wang and E. Chen. Variational principles for BS dimension of subsets. \emph{Dyn. Syst.} 27 (2012), no. 3, 359–385.

\bibitem[23]{re23} T. Wang. Some notes on topological and measure-theoretic entropy. \emph{Qual. Theory Dyn. Syst.} 20 (2021), no. 1, Paper No. 13, 13 pp.

\bibitem[24]{re30} Y. Wang and Zhao. C. Localized topological pressure for amenable group actions. \emph{Anal. Math. Phys.}12(2022), no.3, Paper No. 74, 14 pp.

\bibitem[25]{re3} R. Zhang. Topological pressure of generic points for amenable group actions. \emph{J. Dynam. Differential Equations}30(2018), no.4, 1583–1606.

\bibitem[26]{re6} C. Zhao, E. Chen, X. Hong and X. Zhou. A formula of packing pressure of a factor map. \emph{Entropy}19(2017), no.10, Paper No. 526, 9 pp.

\bibitem[27]{re5} D. Zheng and E. Chen. Bowen entropy for actions of amenable groups. \emph{Israel J. Math.} 212 (2016), no. 2, 895–911.

\bibitem[28]{re4} D. Zheng and E. Chen. Topological entropy of sets of generic points for actions of amenable groups. \emph{Sci. China Math.} 61 (2018), no. 5, 869–880.

\bibitem[29]{re2} X. F. Zhong and Z. J. Chen. Variational principles for topological pressures on subsets. \emph{Nonlinearity} 36 (2023), no. 2, 1168–1191.

\bibitem[30]{re7} X. Zhou, E. Chen,W. C. Cheng. Packing entropy and divergence points. \emph{Dyn. Syst.} 27 (2012), no. 3, 387–402.

\end{thebibliography}
\end{document}